%% file: conebend.tex
\documentclass[a4paper]{article}

\usepackage{amsfonts,amssymb,amsthm}
\usepackage[totalwidth=17cm,totalheight=22cm]{geometry}
\usepackage{psfrag}
\usepackage[T1]{fontenc}

\newcommand{\N}{{\mathbb N}}
\newcommand{\R}{{\mathbb R}}
\newcommand{\Z}{{\mathbb Z}}

\newcommand{\cH}{{\mathcal H}}
\newcommand{\cM}{{\mathcal M}}
\newcommand{\cT}{{\mathcal T}}
\newcommand{\cML}{{\mathcal {ML}}}
\newcommand{\QF}{{\mathcal {QF}}}
\newcommand{\Mt}{\tilde{M}}

\newcommand{\Tt}{\tilde{T}}
\newcommand{\gt}{\tilde{g}}
\newcommand{\gb}{\overline{g}}

\newcommand{\Sb}{\overline{S}}
\newcommand{\taub}{\overline{\tau}}

\newcommand{\be}{\begin{equation}}
\newcommand{\ee}{\end{equation}}

\newtheorem{prop}{Proposition}[section]
\newtheorem{df}[prop]{Definition}
\newtheorem{lemma}[prop]{Lemma}

\newtheorem{thm}[prop]{Theorem}
\newtheorem{cor}[prop]{Corollary}

\newtheorem{remark}[prop]{Remark}

\newtheorem{question}[prop]{Question}
\newtheorem{claim}[prop]{Claim}

\newcommand{\dr}{\partial}

\newcommand{\isom}{\mbox{Isom}}

\newcommand{\sigmat}{\tilde{\sigma}}

\begin{document}

\title{The convex core of quasifuchsian manifolds with particles}

\date{February 2013 (v2)}

\author{
Cyril Lecuire
\thanks{
Institut de Math\'ematiques de Toulouse, UMR CNRS 5219,
Universit{\'e} Toulouse III,
31062 Toulouse Cedex 9,
France. Partially supported by the ANR programs Groupes and 
ETTT (ANR-09-BLAN-0116-01).}
and
Jean-Marc Schlenker
\thanks{
Institut de Math\'ematiques de Toulouse, UMR CNRS 5219,
Universit{\'e} Toulouse III,
31062 Toulouse Cedex 9,
France. Partially supported by the ANR program ETTT (ANR-09-BLAN-0116-01).}
}

\maketitle

\begin{abstract}

We consider quasifuchsian manifolds with ``particles'', i.e., cone
singularities of fixed angle less than $\pi$  going from one connected component
of the boundary at infinity to the other. Each connected component of the
boundary at infinity is then endowed with a conformal structure marked by the
endpoints of the particles. We prove that this defines a homeomorphism from 
the space of quasifuchsian metrics with $n$ particles (of fixed angle) and the
product of two copies of the Teichm\"uller space of a surface with $n$ marked
points. This extends the Bers Double Uniformization theorem to quasifuchsian
manifolds with ``particles''.

Quasifuchsian manifolds with particles also have a convex core. Its boundary
has a hyperbolic induced metric, with cone singularities at the intersection with
the particles, and is pleated along a measured geodesic lamination. 
We prove that any two hyperbolic metrics with cone
singularities (of prescribed angle) can be obtained, and also that any two 
measured bending laminations, satisfying some obviously necessary conditions,
can be obtained, as in \cite{bonahon-otal} in the non-singular case. 

\bigskip

\begin{center} {\bf R{\'e}sum{\'e}} \end{center}

On considère des variétés quasifuchsiennes ``à particules'', c'est-à-dire ayant des
singularités coniques d'angle fixé inférieur à $\pi$ allant d'une composante
connexe à l'infini à l'autre. Chaque composante connexe du bord à l'infini
est alors muni d'une structure conforme marquée par les extrémités des
particules. On montre que ceci définit un homéomorphisme de l'espace
des métriques quasifuchsiennes à $n$ particules (d'angle fixé) vers le produit
de deux copies de l'espace de Teichmüller d'une surface à $n$ points marqués. 
Ceci étend le théorème de double uniformisation de Bers aux variétés quasifuchsiennes
à ``particules''.

Les variétés quasifuchsiennes à particules ont aussi un coeur convexe. Son bord
a une métrique induite hyperbolique, avec des singularités coniques aux 
intersections avec les particules, et est plissé le long d'une lamination
géodésique mesurée. On montre que toute paire de métriques
hyperboliques à singularités coniques (d'angle prescrit) peut être obtenu,
et aussi que toute paire de laminations de plissages, satisfaisant des
conditions clairement nécessaires, peut être obtenus, comme dans le
cas non-singulier \cite{bonahon-otal}.

\end{abstract}

\maketitle

\vspace{0.4cm}

%{\bf AMS classifications:} 

\vspace{0.4cm}

%{\bf Key-words:} 

\bigskip

\tableofcontents

\input conebend1
% intro

\input conebend2
% geometry of the convex core

\input conebend3
% compactness / bending laminations

\input conebend4
% prescribing the bending lamination

\input conebend5

% quasi-conformal estimates

\input conebend6
% the conformal structure at infinity 
% including induced metric 

\input conebend7
% remarks etc

\appendix

\input conebend8

\section*{Acknowledgements}

\bibliographystyle{alpha}
\bibliography{../../papiers/outils/biblio}

\end{document}

%% file: conebend1.tex
\section{Introduction}

\subsection{Convex co-compact manifolds with particles}

\paragraph{Quasifuchsian manifolds.}

A quasifuchsian manifold is a complete hyperbolic manifold $M$ diffeomorphic
to $S\times \R$, where $S$ is a closed, oriented surface of genus at least $2$,
which contains a non-empty, compact, geodesically convex subset, see \cite{thurston-notes}.
Such a manifold has a boundary at infinity, which is the union of two copies
of $S$. Each of those two copies has a conformal structure, $\tau_+$ and $\tau_-$,
induced by the hyperbolic metric on $M$. 
A celebrated theorem of Bers \cite{bers,ahlfors-bers} asserts that 
the map sending a quasifuchsian metric to $(\tau_+,\tau_-)$ 
determines a parameterization of the space of quasifuchsian metrics
on a $M$ by the product of two copies of the Teichm\"uller space of $S$.

A quasifuchsian manifold $M$ contains a smallest non-empty geodesically convex subset,
called its convex core $C(M)$. Here we say that a subset $K\subset M$ is geodesically
convex if any geodesic segment in $M$ with endpoints in $K$ is contained in $K$. This
implies that the inclusion of $K$ in $M$ is a homotopy equivalence. 

The boundary of $C(M)$ is again the union of two copies
of $S$, and each is a pleated surface in $M$, with a hyperbolic induced metric $m_+,m_-$ 
and a measured bending lamination $\lambda_+,\lambda_-$. 
(There is a special, ``Fuchsian'' case, where $C(M)$ is a totally geodesic surface,
the two pleated surfaces mentioned here are then the same, $m_+=m_-$, and 
$\lambda_+=\lambda_-=0$.) 

It is known that any two hyperbolic
metrics can be obtained in this way, this follows from \cite{epstein-marden} or from 
\cite{L4}, however it is not known whether any couple $(m_+,m_-)$ can be uniquely
obtained. 
Similarly, any two measured laminations $\lambda_+,\lambda_-$ can be obtained in this
manner if the weight of any leaf is less than $\pi$ and if $\lambda_+$ and $\lambda_-$
fill $S$ \cite{bonahon-otal}, but it is not known whether uniqueness holds. Recall
that $\lambda_-$ and $\lambda_+$ fill $S$ if there exists $\epsilon>0$ such that, for 
any closed curve $c$ in $S$, $i(\lambda_-,c)+i(\lambda_+,c)\geq \epsilon$.

The main goal here is to extend those results to quasifuchsian manifols with ``particles'',
that is, cone singularities of a certain type connecting the two connected components of
the boundary at infinity, as described below.

Note that all the results mentioned here are actually known in the more general
context of convex co-compact hyperbolic manifolds, i.e., interiors of compact manifolds with
boundary, with a complete hyperbolic metric, containing a non-empty, compact, geodesically
convex subset (or, even more generally for geometrically finite hyperbolic manifolds), the
result concerning the measured bending lamination of the boundary can then be found in 
\cite{lecuire}. We stick here to the quasifuchsian setting for simplicity.

\paragraph{Cone-manifolds.}

We consider here hyperbolic cone-manifolds of a special kind, which have cone
singularities along curves (a more general notion is defined in \cite{thurston-notes},
allowing for singularities along graphs). Let $\theta\in (0,\pi)$, we call 
$H^3_\theta$ the hyperbolic manifold with cone singularities obtained by gluing
isometrically the two faces of a hyperbolic wedge of angle $\theta$ (the closed
domain in $H^3$ between two half-planes having the same boundary line). There is
a unique such gluing which is the identity on the ``axis'' of the wedge. We will
be using here the following (restrictive) definition.

\begin{df}
A hyperbolic cone-manifold is a manifold along with a metric for which each point 
has a neighborhood modeled on $H^3_\theta$ for some $\theta\in (0,\pi)$. 
\end{df}

Let $M$ be a hyperbolic cone-manifold, it has two kind of points. Those which
have a neighborhood isometric to a neighborhood of a point of some $H^3_\theta$
outside the cone singularity are called regular points, while the others are
called singular points or cone points. The set of regular points will be denoted
by $M_r$, and the set of singular points by $M_s$. By definition, $M_s$ is a
union of curves, if $M$ is complete then those curves can be either closed 
curves or infinite lines. To each of those curves is associated an angle 
$\theta\in (0,\pi)$ --- such that all points have a neighborhood isometric to
a neighborhood of the cone singularity in $H^3_\theta$ --- which is called its
cone angle or simply its angle.

Recall the usual notion of convexity, which differs from other possible
notions (e.g. the local convexity of the boundary of a domain). 

\begin{df}
Let $M$ be a hyperbolic cone-manifold. A subset $K\subset M$ is {\bf 
geodesically convex}
if any locally geodesic segment in $M$ with endpoints in $K$ is contained in $K$.
\end{df}

A non-empty geodesically convex subset of $M$ is homotopically equivalent to $M$
and contains all closed geodesics of $M$, see \cite[Lemma A.12]{qfmp}. 

\paragraph{Quasifuchsian manifolds with particles.}

Quasifuchsian manifolds with particles are defined in the same way as non-singular
quasifuchsian manifold. 

\begin{df}
A {\bf quasifuchsian manifold with particles} is a complete
hyperbolic cone-manifold $M$ isometric to the product $S\times \R$, 
where $S$ is a closed, orientable surface, endowed
with a complete hyperbolic metric with cone singularities of angles in $(0,\pi)$
on the lines $\{ x_i\}\times \R$, for $x_1,\cdots, x_{n_0}$ distinct points in $S$,
which contains a non-empty, compact, geodesically convex subset. 
We require that $n_0\geq 4$ if $S$ is a sphere, i.e. that $M$ has at least $4$ singularities,
and that $n_0\geq 1$ if $S$ is a torus. %% j'ai ajoute cette condition pour le tore, on 
%% pourrait en principe prendre n_0=0 pour les autres cas et autant l'inclure donc.
\end{df}

Notice that the definition also makes sense if $S$ is a sphere and $n_0=3$ but then 
the metric would be uniquely defined (up to isotopy) by the cone angles. 
We are not considering this case for technical reasons but also because there is not much to say about it.

Given a non-empty convex subset $K$ of a quasifuchsian manifold with particles,
then $K$ contains all closed geodesics of $M$ (see \cite[Lemma A.12]{qfmp}) and
the inclusion of $K$ in $M$ is a homotopy equivalence (this is proved below). 

Geometrically, quasifuchsian manifolds with particles can be considered as $I$-bundles
in the category of hyperbolic 3-manifolds with cone singularities. The term ``particle''
comes from physical motivations. Quasifuchsian manifolds have Lorentzian siblings
called Anti-de Sitter (AdS) globally hyperbolic (GH) manifolds which share many of
the key properties recalled above, see \cite{mess,mess-notes}. From a physics viewpoint,
GH AdS 3-manifolds are a 3-dimensional toy model for gravity, as they model
an empty space with negative cosmological constant. To go beyond an empty
model, massive point particles can be added and modeled as cone singularities
along time-like lines, see e.g. \cite{thooft1,thooft2}. The resulting GH AdS 
manifolds with particles display some properties which are parallel to those
obtained here, see \cite{cone}.

The restrictions on the cone angles --- supposed to be in $(0,\pi)$ --- are necessary
at several points here, as they were in \cite{qfmp}. It seems to be physically relevant,
too. We will mention some points where this hypothesis is useful below as they occur.
We do not know whether Theorem \ref{tm:bers}, for instance, can be extended to 
cone angles less than $2\pi$. In the parallel Lorentzian theory concerning globally hyperbolic
anti-de Sitter manifolds, new phenomena arise when the cone angles are larger than $\pi$, 
see \cite{colI,colII}.

Quasifuchsian manifolds with particles are always considered here up to 
isotopies.

\paragraph{Convex co-compact manifolds with particles.}

The previous definition can be extended to a definition of convex co-compact 
manifolds with particles.

\begin{df}
A {\bf convex co-compact hyperbolic manifold with particles} is a complete
hyperbolic cone-manifold $M$ such that:
\begin{itemize}
\item $M$ is homeomorphic to the interior of a compact
manifold with boundary $N$, with non-trivial fundamental group,
\item the singular locus corresponds under the homeomorphism to a disjoint 
union of curves in $N$ with endpoints on $\dr N$,
\item the angle at each singular curve is less than $\pi$,
\item $M$ contains a non-empty compact subset which is convex.
\end{itemize}
\end{df}

A further extension to geometrically finite manifolds with particles
is possible, we leave the details to the interested reader. We consider
here only quasifuchsian mainfolds (with particles) although some of the
intermediate statements can be extended to convex co-compact manifolds
with particles. There is also some hope to extend the main results to
this more general setting, however some technical hurdles have to be
overcome before this can be achieved.

%%paragraphe ajoute
Note that there is another possible notion of quasifuchsian manifolds
with cone singularities: those which are singular along closed curves,
as studied in particular by Bromberg \cite{bromberg2,bromberg1}. Although
there are similarities between those two kinds cone-manifolds (in particular
concerning their rigidity), the questions considered here are quite different
from those usually associated to those considered for quasifuchsian
cone-manifolds with singularities along closed curves (drilling of
geodesics, etc).

\subsection{The conformal structure at infinity.}

\paragraph{Conformal structures and hyperbolic metrics on surfaces.}

Let's fix some notations.

\begin{df}
Let $S$ be a closed 
orientable %% ajoute
surface, let $x_1, \cdots, x_{n_0}\in S$ be distinct points with $n_0\geq 4$ if $S$ is a sphere 
and $n_0\geq 1$ is $S$ is a torus, and let $\theta_1, \cdots, \theta_{n_0}\in (0,\pi)$.
We then call:
\begin{itemize}
\item $\cT_{S,x}$ the space of conformal structures on $S$, considered up
to isotopies of $S$ fixing the $x_i$,
\item $\cH_{S,x,\theta}$ the space of hyperbolic metrics on $S$, with cone
singularities at the $x_i$ where the angle is $\theta_i$, considered
up to isotopies fixing the $x_i$.
\end{itemize}
\end{df}

There is a one-to-one map between $\cT_{S,x}$ and $\cH_{S,x,\theta}$, 
because any conformal structure contains a unique hyperbolic metric
with cone singularities at the $x_i$  of prescribed angle (see 
\cite{troyanov}). We keep distinct notations for clarity.

Notice again that these definition also make sense when $S$ is a sphere and $n_0=3$ 
provided that $\sum_{i=1}^{n_0} \theta_i -2\pi <0~$. In this case the spaces $\cT_{S,x}$ and $\cH_{S,x,\theta}$ are points. 

The statements considered here are already well understood when $n_0=0$, so we will
focus below on the case $n_0\geq 1$.

\paragraph{The conformal structure at infinity.}

Non-singular quasifuchsian manifolds have a natural conformal structure
at infinity, which can be defined by considering the action of their fundamental
group on their discontinuity domain, see \cite{thurston-notes}. This definition
cannot be used directly for quasifuchsian manifolds with particles, however it
is still possible to define a conformal structure at infinity, see \cite[Section 3.2]{qfmp}.

Therefore, to each quasifuchsian metric $g\in \QF_{\theta_1,\cdots,\theta_{n_0}}$ are
associated two points $\tau_+,\tau_-\in \cT_{S,n_0}$ corresponding
to the conformal structures --- marked by the endpoints of the ``particles'' ---
on the upper, resp. lower, connected component of the boundary at infinity.

Note that we always implicitly consider conformal structures on the boundary
at infinity up to isotopy. (It is therefore not necessary to consider markings.)

\paragraph{A compactness lemma for the conformal structure at infinity.}

We consider again a closed surface $S$ along with $n_0$ distinct points ($n_0\geq 1$) 
$x_1,\cdots, x_{n_0}\in S$
and angles $\theta_1,\cdots, \theta_{n_0}\in (0,\pi)$ so that
$$ 2\pi \chi(S) - \sum_{i=1}^{n_0} (2\pi-\theta_i) <0~. $$

\begin{prop} \label{pr:compact-conforme}
Let $(g_n)_{n\in \N}$ be a sequence of quasifuchsian metrics on 
$S\times \R$, with particles (cone singularities) on the lines 
$\{ x_i\}\times \R$, of angle equal to $\theta_i$. Suppose that
the conformal structures at infinity, $\tau_{-,n}, \tau_{+,n}\in \cT_{S,n_0}$, 
converge to conformal structures $\tau_{-,\infty}, \tau_{+,\infty}$. 
Then $(g_n)_{n\in \N}$ has a subsequence converging to a quasifuchsian
metric with particles.
\end{prop}

The proof is contained in Section \ref{ssc:62}, it is based on 
the compactness results described below (in Section \ref{se:compact-lamination}) 
relative to the induced metric and bending lamination on the boundary of the convex core. 

\paragraph{A Bers-type theorem with particles.}

Using the previous proposition, along with the main result of \cite{qfmp},
leads to an extension to quasifuchsian manifolds with particles of
a classical result of Bers \cite{bers} on ``double uniformization''. 

\begin{thm} \label{tm:bers}
The map from $\QF_{S,x,\theta}$ to $\cT_{S,x}\times \cT_{S,x}$
sending a quasifuchsian hyperbolic metric to the conformal structures
at $+\infty$ and at $-\infty$ (marked by the endpoints of the particle) is a homeomorphism.
\end{thm}

\subsection{The geometry of the convex core}

\paragraph{Measured laminations.}

We refer the reader to e.g. \cite{Casson,penner-harer,otal-hyperbolisation} 
for the definition and main properties
of measured laminations on closed (non-singular) surfaces as well as the 
topology on the space of measured laminations. There
are two possible definitions. One is geometric, in terms of measured geodesic
laminations on hyperbolic surfaces, with a transverse measure, while the
other definition is topological, and can involve the boundary at infinity of
the universal cover of %% ajoute
the surface. The two definitions are equivalent, basically because, in 
a closed (or finite volume) hyperbolic surface, any closed curve can be
realized uniquely as a closed geodesic. 

\begin{prop}
Let $\Sigma$ be a hyperbolic surface with cone singularities, where the angle
is less than $\pi$. Let $\lambda$ be a (topological) lamination on
$\Sigma$. Then $\lambda$ can be realized uniquely as a geodesic lamination.  
\end{prop}

The space of measured geodesic laminations on a hyperbolic surface with cone
singularities of angles less than $\pi$ therefore does not depend on the
cone angles. 

\begin{df}
We call $\cML_{S,n_0}$ the space of measured lamination on $S$ with $n_0$ marked points.
\end{df}

Thus, for any hyperbolic metric $m$ on $S$ with $n$ cone singularities of angle less
than $\pi$,  $\cML_{S,n_0}$ can be canonically identified with the space of measured
geodesic laminations on $(S,m)$.

\paragraph{The convex core.}

The following basic proposition can be found in the appendix of \cite{qfmp}.

\begin{prop}
Let $M$ be a convex co-compact hyperbolic manifold with particles, and let
$K$ and $K'$ be two non-empty geodesically convex subsets. Then $K\cap K'$ is a non-empty
geodesically convex subset.  
\end{prop}

It leads to a natural definition.

\begin{df}
Let $M$ be a quasifuchsian manifold with particles. Its 
{\bf convex core} $C(M)$ 
is the smallest non-empty geodesically convex subset contained in it.
\end{df}

By construction, $C(M)$ is a ``minimal'' convex subset of $M$ and it
follows from general arguments (see \cite{thurston-notes}) that its
boundary is, outside the singular curves, a pleated surface (a locally
convex, ruled surface). It turns out that, under the condition that
the cone angles are less than $\pi$, the boundary of $C(M)$ 
intersects the cone singularities orthogonally, and is even
totally geodesic in the neighborhood of each such intersection, see
\cite[Lemma A.15]{qfmp}. 

It follows that there is a well-defined notion of closest-point projection
from $M$ to $C(M)$. As a consequence, the inclusion of $C(M)$ in $M$
is a homotopy equivalence. The same holds for any non-empty convex 
subset of $M$.

Therefore, given a quasifuchsian metric with particles $g\in \QF_{S,\theta}$, 
the induced metrics on the upper and lower boundary
components of $C(M)$ (which might coincide in special cases)
are two hyperbolic metrics $m_+,m_-\in \cH_{S,x,\theta}$. 
Moreover, those two boundary components are pleated
along measured bending laminations $l_+,l_-\in \cML_{S,n_0}$.

\paragraph{A remark on the hypothesis}
A well-known fact concerning hyperbolic surfaces with cone singularities
is that, as long as the cone angles are less than $\pi$, it remains true
that any homotopy class of closed curves in the regular part contains a unique geodesic
(see e.g. \cite{dryden-parlier}). A fairly direct consequence is that,
as for closed surfaces, any topological measured lamination (in the complement
of the cone singularities) can be uniquely realized as a measured geodesic lamination.

This is one reason --- albeit not the only one --- why it is relevant to consider
here cone singularities of angle less than $\pi$, rather than less than $2\pi$. Indeed
for cone singularities of angles less than $2\pi$, the induced metric on the boundary
of the convex core might also have cone singularities of angle between $\pi$ and 
$2\pi$, and for those metrics the one-to-one relation between measured laminations
and measured geodesic laminations is lost.

\subsection{Prescribing the bending lamination}

\paragraph{Results in the non-singular case.}

For non-singular convex co-compact hyperbolic manifolds an existence and
uniqueness theorem for metrics with a given rational measured bending lamination
was proved by Bonahon and Otal \cite{bonahon-otal}. 
(Recall that a lamination is rational if its support is a disjoint union of
closed curves.)
When the lamination is 
not rational, an existence result was proved in \cite{bonahon-otal} for manifolds
with incompressible boundary, it was extended in \cite{lecuire} to manifolds
with compressible boundary. 

\paragraph{Rational laminations with particles.}

As for quasifuchsian manifolds (without particles), 
it is possible to give an existence and uniqueness statement concerning
the bending lamination on the boundary of the convex core, but only
for rational laminations.

\begin{thm} \label{tm:rational}
Let $S$ be a closed 
orientable %% ajoute
surface, let $x_1, \cdots, x_{n_0}\in S$ be distinct
points, and let $\theta_1, \cdots, \theta_{n_0}$ be in $(0,\pi)$. 
Suppose that $n_0\geq 4$ if $S$ is a sphere, and that $n_0\geq 1$ if $S$ is a torus. %% ajoute
Let $\lambda_-,\lambda_+
\in \cML_{S,x}$ be measured laminations, each with support a disjoint union
of closed curves. Suppose that: 
\begin{itemize}
\item $\lambda_-$ and $\lambda_+$ fill $S$,
\item each closed curve in the support of $\lambda_-$ (resp. $\lambda_+$)
has weight less than $\pi$.
\end{itemize}
Then there exists a metric $g\in \QF_{S,x,\theta}$ such
that the measured bending lamination on the upper (resp. lower)
boundary component of the convex core of $(S\times \R, g)$ is 
$\lambda_+$ (resp. $\lambda_-$). Moreover $g$ is unique up to isotopies.
\end{thm}

The proof, which is given in Section \ref{se:rational}, is based
on the rigidity theorem of Hodgson and Kerckhoff \cite{HK} for
closed hyperbolic manifolds with cone singularities. We prove in
Lemma \ref{lm:n-general} that the hypothesis are necessary conditions.

\paragraph{General laminations.}

When considering laminations which are not necessarily rational, we
obtain only a weaker result, because we can only claim existence,
but not uniqueness (this remains an open problem even in the non-singular
case, see \cite{bonahon-otal,lecuire}).

\begin{thm} \label{tm:general}
Let $S$ be a closed surface, let $x_1, \cdots, x_{n_0}\in S$ be distinct
points, and let $\theta_1, \cdots, \theta_{n_0}$ be in $(0,\pi)$. 
Let $\lambda_-,\lambda_+ \in \cML_{S,x}$. Suppose that: 
\begin{itemize}
\item $\lambda_-$ and $\lambda_+$ fill $S$,
\item each closed curve in the support of $\lambda_-$ (resp. $\lambda_+$)
has weight less than $\pi$.
\end{itemize}
Then there exists a metric $g\in \QF_{S,x,\theta}$ such
that the measured pleating lamination on the upper (resp. lower)
boundary component of the convex core of $(S\times \R, g)$ is 
$\lambda_+$ (resp. $\lambda_-$).
\end{thm}

The two conditions on $\lambda_-, \lambda_+$ in this theorem are easily seen to be necessary when $g$ is not fuchsian,
see Lemma \ref{lm:n-general}. Note that both Theorem \ref{tm:rational} and
Theorem \ref{tm:general} are restricted to quasifuchsian manifolds with particles,
rather than more general convex co-compact manifolds with particles.

\subsection{The induced metric on the boundary of the convex core.}

The Bers-type result on the conformal metric at infinity can
be used to obtain an existence result concerning
the prescription of the induced metric on the boundary of the 
convex core. 

\begin{thm} \label{tm:metriques}
Let $m_-, m_+\in \cH_{S,n_0,\theta}$, where $\theta=(\theta_1,\cdots,\theta_{n_0})\in 
(0,\pi)^n$. There exists a quasifuchsian 
metric with particles on $S\times \R$, with particles of angle
$\theta_i$ at the lines $\{ x_i\}\times \R$, for which the induced
metric on the boundary of the two connected components of the convex 
core are $m_-$ and $m_+$.
\end{thm}

In the smooth case -- i.e. for quasifuchsian hyperbolic manifolds without
conical singularities -- the corresponding result is well-known, it follows 
either from results of Labourie \cite{L4} or from a
partial answer, first given by Epstein and Marden \cite{epstein-marden}, to a
conjecture of Sullivan. (The conjecture made by Sullivan turned out to be 
wrong, see \cite{epstein-markovic}, but the result proved by Epstein and Marden is sufficient
to prove Theorem \ref{tm:metriques} in the non-singular context.)

As for the conformal structure at infinity, it might be possible to
extend this statement to cover convex co-compact (resp. geometrically
finite) manifolds with particles. 
The uniqueness remains elusive, as in the non-singular case. 

\subsection{Applications}

Quasifuchsian manifolds can be used as tools in Teichm\"uller theory. By
extension, the quasifuchsian manifolds with particles considered here 
can be used as tools for the study of the Teichm\"uller space of 
hyperbolic metrics with cone singularities (of angle less than $\pi$)
on a surface. 

One such application is through the renormalized volume of those
quasifuchsian manifolds with particles, as considered in \cite{volume,review}.
In the non-singular case this renormalized volume is equal to the 
Liouville functional (see \cite{TZ-schottky,takhtajan-teo,takhtajan-zograf:spheres}), 
it is a K\"ahler potential on $\cH_{S,x,\theta}$. Other applications of closely
related tools, in the non-singular context, for the global geometry of the Weil-Petersson
metric on Teichm\"uller space, can be found in \cite{McMullen}.
Yet other applications, to some properties of the grafting map, are considered
in \cite{cp}, and the manifolds with particles considered here should allow for
an extension to the grafting map on $\cH_{S,x,\theta}$. 

\subsection{Outline of the proofs}

We now turn to a description of the main technical points of the proofs. 

\paragraph{Measured bending laminations.}

Theorem \ref{tm:rational} is proved by an argument strongly influenced
by the proof given by Bonahon and Otal \cite{bonahon-otal} for non-singular
convex co-compact manifolds. Thanks to a doubling trick, the infinitesimal
rigidity of the convex cores of convex co-compact manifolds with particles, 
with respect to the (rational) measured bending lamination, is reduced to
an important infinitesimal rigidity result proved for hyperbolic cone-manifolds
by Hodgson and Kerckhoff \cite{HK}. A deformation argument then provides the
proof of the theorem. 

The existence result for general laminations on quasifuchsian manifolds with
particles (Theorem \ref{tm:general}) 
can then be obtained by an approximation argument, as in the 
non-singular case in \cite{bonahon-otal,lecuire}. The key step of the proof
is a compactness statement, showing that if the measured bending
laminations converge to a limit having good properties, then the quasifuchsian
metrics converge after extracting a subsequence. However the arguments 
developed in \cite{bonahon-otal,lecuire} cannot be used in the context
of quasifuchsian manifolds with particles, because they rely heavily on 
the representation of the fundamental group. Different arguments are
therefore used here, which are more differential-geometric in nature. 

Those arguments are sometimes technically involved because
of the added difficulties induced by the particles. However, after stripping
the proof of the elements which are needed only because of the particles
(for instance the multiple cover argument used in Section \ref{ssc:34}
to find simplicial surfaces with given boundary in the convex core), 
the compactness proof given here is simpler than the one in \cite{bonahon-otal,lecuire}.

\paragraph{Prescribing the induced metric on the boundary of the convex core.}

We give in Section \ref{se:induced} a rather elementary proof 
of Theorem \ref{tm:metriques}, which has two parts. Call 
$t_-$ (resp. $t_+$) the hyperbolic metric in the conformal class
$\tau_-$ (resp. $\tau_+$) with cone angles $\theta_i$ at
$x_i$. The first part 
is an upper bound on the length of the curves in the hyperbolic
metric at infinity $t_\pm$, following \cite{bridgeman,bridgeman-canary}.

Now recall that, by Thurston's Earthquake Theorem 
\cite{kerckhoff,thurston-earthquakes}, there exists a unique 
right earthquake sending a given hyperbolic metric to another one.
This extends to hyperbolic metrics with cone singularities of angle
less than $\pi$, see \cite{cone}. In particular there is a unique measured
lamination $\nu_+$ such that the right earthquake along $\nu_+$,
applied to $m_+$, yields the hyperbolic metric $t_+$. 
The second part of our proof is a bound on the length
of $\nu_+$ for $m_+$ (see Proposition \ref{pr:length}).

This is then used in Section \ref{se:induced}
to prove Theorems \ref{tm:bers} and \ref{tm:metriques}. The proof of
Theorem \ref{tm:bers} also uses another main ingredient, the local rigidity
of quasifuchsian manifolds with particles proved in \cite{qfmp}.

\paragraph{Quasi-conformal estimates.}

There is another possible way to prove Theorem \ref{tm:metriques},
closer to the argument used in the non-singular case (as seen in
\cite{epstein-marden,bridgeman,bridgeman-canary}). It uses a 
bound on the quasi-conformal factor between the 
conformal structure at infinity $\tau$
and the conformal class of the induced metric $m$ on the boundary of 
the convex core, both understood as elements of $\cT_{S,n_0}$. 

\begin{prop} \label{pr:qconf}
There exists a constant $C>0$ (depending only on the topology
of $M$) such that $\tau$ is $C$-quasiconformal to $m$.
\end{prop}

This proposition is not formally necessary to obtain the main results
presented here, its proof can be found in Appendix \ref{se:quasiconf}.

As mentioned above, the proof of Theorem \ref{tm:metriques} through 
Proposition \ref{pr:qconf} would be much closer to the proof(s) known
in the non-singular case. It can be pointed out that the proof given
in Section \ref{se:induced} is quite parallel, but in the context of
Teichm\"uller theory understood as the study of hyperbolic rather than 
complex surfaces. From this viewpoint, Proposition \ref{pr:length}
is a direct analog of Proposition \ref{pr:qconf}, with quasiconformal
deformations replaced by earthquakes.

\paragraph{What follows.}

Section 2 presents the definition of the convex core of a convex co-compact
manifold with particles, and some of its simple properties, extending
well-known properties with no cone singularity. In Section \ref{se:compact-lamination} we state and
prove a key compactness statement with respect to the measured bending 
lamination on the boundary of the convex core. 
Section \ref{se:rational} contains the proof of Theorem
\ref{tm:rational}, using a local rigidity statement of Hodgson and Kerckhoff
\cite{HK} and the compactness Lemma of Section \ref{se:compact-lamination}. 
Section \ref{se:earthquake} contains
the proof of Theorem \ref{tm:general}, and Section \ref{se:induced} contains the proof
of Theorem \ref{tm:bers} and of Theorem \ref{tm:metriques}. Section \ref{se:questions}
contains some remarks on the analogy with corresponding problems in 
anti-de Sitter geometry and on applications to the Weil-Petersson metric
of the Teichm\"uller space of hyperbolic metrics with cone singularities
of prescribed angles on a closed surface (see \cite{volume,review,cp}).
Finally, Appendix \ref{se:quasiconf} contains the proof of Proposition
\ref{pr:qconf}, based on the estimates on the length of the earthquake
lamination obtained in Section \ref{se:earthquake}.

%% file: conebend2.tex
\section{The geometry of the convex core}
\label{se:convex-core}

This section contains some basic statements necessary to understand 
the geometry of convex co-compact manifolds with particles, concerning
in particular the convex core and its boundary. We consider here
such a convex co-compact manifold with particles, $M$, and denote by
$M_r$ its regular part and by $M_s$ its singular part (the union of the
singular lines). 

We exclude below the simplest case where $M$ is Fuchsian, that is, where
it is the warped product of a hyperbolic surfaces with cone singularities
$(S,h)$ by $\R$, with the metric $dt^2+\cosh(t)^2 h$. In this Fuchsian
case the convex core is a surfaces, corresponding to $t=0$, and it
is totally geodesic outside the intersection with the particles, and 
orthogonal to those particles.

\subsection{Surfaces orthogonal to the singular locus}

We define here a natural notion of pleated surface orthogonal to the 
singular locus in $M$. The first step is to define the
notion of totally geodesic plane orthogonal to a cone singularity
in a hyperbolic cone-manifold. The first condition is that the
surface is totally geodesic outside its intersections with the 
particles. The second condition is local, in the neighborhood of  
the intersections with the particles; there, the surface should
correspond to the image in $H^3_\theta$ of the restriction to 
the wedge (used to define $H^3_\theta$) of a plane orthogonal
to the axis of the wedge.

\begin{df}
Let $\Sigma$ be a pleated surface in $M_r$, and let $\Sigma'$ be its closure as a
subset of $M$; suppose that $\Sigma'\setminus \Sigma\subset M_s$. 
We say that $\Sigma'$ is {\it orthogonal to the singular locus}
if any $x\in \Sigma'\setminus \Sigma$ has a neighborhood in $\Sigma'$
which is a totally geodesic surface orthogonal to the singular locus. 
\end{df}

This definition can be extended to encompass more general surfaces, i.e., 
surfaces which are neither pleated nor totally geodesic in the neighborhood
of the singular locus. In this more general case the definition can be 
given in terms of the convergence of the unit normal vector to a
vector ``tangent'' to the singular locus at its intersection with the
surface. This will however not be needed here.

\subsection{The convex core of a manifold with particles}

Among the defining properties of a quasifuchsian cone-manifold
$M$ is the fact that it contains a compact subset $K$ which is convex 
in the (strong) sense that any geodesic segment in $M$ with endpoints
in $K$ is contained in $K$. We have already seen that it is possible to
define the {\it convex core} of $M$ as the smallest compact subset of 
$M$ which is convex, denoted by $C(M)$.

\begin{thm} \label{tm:cc}
Suppose that $C(M)$ is not a totally geodesic surface. Then its boundary
is the disjoint union of surfaces which are orthogonal to the 
singular locus. Each connected component of the singular locus of $M$
intersects $C(M)$ along a segment.
\end{thm}

The proof is a consequence of two lemmas, both stated under the hypothesis of
the theorem. The second lemma in particular gives more precise informations on
the geometry of the convex core, it is taken from \cite[Lemma A.14]{qfmp}.

\begin{lemma} \label{lm:cc-ortho}
The boundary of $C(M)$ is a surface orthogonal to the singular locus.
\end{lemma}

Let $x\in M$, we denote by $L_x$ the {\bf link} of $M$ at $x$, that is, the
space of geodesic rays starting from $x$ (parametrized at speed $1$), with
its natural angle distance. When $x$ is a regular point of $M$, $L_x$ is 
isometric to the 2-dimensional sphere $S^2$, with its round metric. When 
$x$ is contained in a singular line of angle $\theta$, $L_x$ can be
described as the metric completion of the quotient by a rotation of angle 
$\theta$ of the universal
cover of the complement of two antipodal points in $S^2$.

\begin{df}
Let $K\subset M$ be convex, and let $x\in M$. The {\bf link} of $K$ at $x$
is the set of vectors $v\in L_x$ such that there is a (small) 
geodesic ray starting from $x$ in the direction of $v$ which is contained in
$K$. It is denoted by $L_x(K)$.
\end{df}

Clearly $L_x(K)=\emptyset$ when $x$ is not contained in $K$, while
$L_x(K)=L_x$ when $x$ is contained in the interior of $K$.

To go further, we define the oriented normal bundle of $\dr C(M)$, denoted by
$N^1_r\dr C(M)$, as the set of $(x,n)\in TM$ such that $x\in C(M)$ is not in
the singular locus of $M$ and $n$ is a unit vector such that its orthogonal is
a support plane of $C(M)$ at $x$, and $n$ is oriented towards the exterior of
$C(M)$. 

Let $x\in M$ be a non-singular point, let $v\in T_xM$ and let $t\in \R_+$. 
For $t$ small enough, it is possible to define the image of $(x,tv)$ by
the exponential map, it is the point $\exp(x,tv):=g(t)$, where $g$ is
the geodesic, parametrized at constant speed, such that $g(0)x$ and 
$g'(0)=v$. As $t$ grows, $\exp(x,tv)$ remains well-defined until $g$
intersects the singular set of $M$. 

\begin{lemma} \label{lm:ends}
The exponential map is a homeomorphism from $N^1_r\dr C(M)\times (0,\infty)$ 
to the complement of $C(M)$ in $M$, and its restriction to the complement of
the points of the form $(x,v,t)$, for $x\in M_s$ and $v$ a singular point of
$L_x$, is a diffeomorphism to complement of $C(M)$ in $M_r$. The map:
$$ 
\begin{array}{clcl}
\exp_\infty: & N^1_r\dr C(M) & \rightarrow & \dr_\infty M \\
& (x,v) & \mapsto & \lim_{t\rightarrow \infty} \exp(x,tv)
\end{array} $$
is a homeomorphism from $N^1_r\dr C(M)$ to the complement in $\dr_\infty M$
of the endpoints of the singular curves in $M$.
\end{lemma}

This follows directly from Lemma A.11 in \cite{qfmp}.

The proof of Theorem \ref{tm:cc} clearly follows from Lemma \ref{lm:cc-ortho} and Lemma
\ref{lm:ends}, since Lemma \ref{lm:ends} shows that the cone singularities cannot re-enter
the convex core after exiting it.

\subsection{The geometry of the boundary}

By construction, $C(M)$ is a minimal convex set in $M$, and it follows
as in the non-singular case (see \cite{thurston-notes}) that its boundary 
is a ``pleated surface'' except at its intersections with the singular
curves. 

\begin{lemma} \label{lm:boundary-cc}
The surface $\dr C(M)$ has an induced metric which is hyperbolic (i.e. it has
constant curvature $-1$) with conical singularities at the intersections of
$\dr C(M)$ with the singular curves of $M$, where the total angle is the same
as the total angle around the corresponding singular curve. It is ``pleated''
along a measured lamination $\lambda$ in the complement of the singular points. 
Moreover the distance between the support of $\lambda$ and the intersection 
of the singular set of $M$ with $\dr C(M)$ is strictly positive. 
\end{lemma}

\begin{proof}
Since $C(M)$ is a minimal convex subset, its boundary is locally convex
and ruled, therefore developable (see \cite{spivak} for the Euclidean analog,
or \cite{thurston-notes}) so that its induced metric is hyperbolic. The fact
that its intersection points are conical singularities, with a total angle
which is the same as the total angle around the corresponding singularities,
is a consequence of the fact that $\dr C(M)$ is orthogonal to the
singularities. 

Similarly, the fact that $\dr C(M)$ is pleated along a measured lamination is
a direct consequence of the fact that it is ruled and locally convex,
i.e. that each point in $\dr C(M)$ is in either a complete hyperbolic
geodesics or a totally geodesic ideal triangle. 
The support of $\lambda$ is a disjoint union of embedded maximal geodesics, 
and it is well-known (see e.g.
\cite{dryden-parlier}) that (under
the hypothesis that the angles at the cone singularities are strictly 
less than $\pi$) embedded geodesics remain at positive distance from the
singular locus. So the distance between the support of $\lambda$ and
the singular locus of $\dr C(M)$ is strictly positive.
\end{proof}

\subsection{The distance between the singular curves}

We state and prove here some elementary statements on the distance between singular
points in the boundary of $C(M)$ and between singular curves in $M$. They will be
useful at several points below. 

\begin{lemma} \label{lm:dist2}
Let $\theta\in (0,\pi)$. There exists $\epsilon>0$ and $\rho>0$, depending on 
$\theta$, such that:
\begin{enumerate}
\item in a complete hyperbolic surface with cone singularities of angle less than $\theta$
(not homeomorphic to a sphere),
two cone singularities are at distance at least $\epsilon$,
\item if $D$ is a closed 2-dimensional geodesic disk of radius $\epsilon$ centered at a 
singular point $x_0$ of cone angle $\theta$, and if $\Omega\subset D$ is
a convex subset whose closure intersects the boundary of $D$, then $\Omega$ contains
all points of $D$ at distance at most $\rho$ from $x_0$.
\end{enumerate}
In particular, it follows from point (1) that no embedded geodesic in $D$ can 
come within distance less than $\rho$ from the cone singularity. 
\end{lemma}

\begin{proof}
The first point is well-known, see e.g. \cite{dryden-parlier}. The interested 
reader can construct an elementary proof based on Dirichlet domains, as in 
3-dimensional manifolds in the proof of the second point, below.

For the second point let $x_1\in \dr D\cap \Omega$, and let 
$\gamma$ be the minimizing geodesic segment from $x_0$ to $x_1$. Since 
$D$ contains no other singular point by the first point, the complement
of $\gamma$ in $D$ is isometric to an angular sector in the disk of radius
$\epsilon$ in $H^2$. This angular sector has three vertices, one corresponding
to $x_0$ and the other two corresponding to $x_1$. Since $\theta<\pi$, it
is convex at the vertex corresponding to $x_0$. Let $s$ be the geodesic
segment joining the two vertices corresponding to $x_1$. Then $\Omega$,
being convex, contains the projection in $D$ of the triangle bounded
by $s$ and by the two geodesic segments in the boundary of $D$ joining
$x_0$ to the two vertices projecting to $x_1$. This proves the statement,
with $\rho$ equal to the distance between $x_0$ and $s$.
\end{proof}

We now turn to a similar lemma, but concerning 3-dimensional manifolds with particles. 

\begin{lemma} \label{lm:dist3}
Let $\theta\in (0,\pi)$. There exists $\epsilon>0$ and $\rho>0$, depending on 
$\theta$, such that, if $M$ is a quasifuchsian manifold with particles of angle less than
$\theta$, then any two particles in $M$ are at distance at least $\epsilon$.
\end{lemma}

\begin{proof}
We reason by contradiction, that is, we fix $\theta\in (0, \pi)$
and for any $n>0$ there is a quasifuchsian manifold $M_n$ with particles of angles less
than $\theta$, with two particles $p_n,p_n'$ at distance less than $1/n$ (and no two
particles strictly closer than $p_n$ and $p'_n$). Let 
$s_n$ be the length-minimizing segment between $p_n$ and $p'_n$, and let $x_n$ be
its midpoint. We call $D_n$ the Dirichlet domain in $M_n$ centered at $x_n$.

We call $(M'_n,x_n)$ the pointed cone-manifold obtained by performing on 
$(M_n, x_n)$ a homothety of ratio $1/L(s_n)$, so that $M'_n$ has constant
curvature $L(s_n)^2\leq 1/n^2$. Let $D'_n$ be the Dirichet domain centered at
$x_n$ in $M'_n$, so that $D'_n$ is obtained by performing a homothety of ratio 
$1/L(s_n)$ on $D_n$.

By construction the cone singularities in $M'_n$ are at distance at least $1$,
so that, after extracting a sub-sequence, $(M'_n, x_n)$ converges in the 
pointed Gromov-Hausdorff topology to a pointed manifold $(M',x)$. Still by
construction, $M'$ contains at least two cone singularities $p$ and $p'$,
limits respectively of $p_n$ and $p'_n$, at distance $1$, with $x$ at the
midpoint of a geodesic segment of minimal length connecting $p$ to $p'$.
Let $D'$ be the Dirichlet domain centered at $x$ in $M'$, then $(D',x)$ 
is the limit of the $(D'_n,x_n)$ in the pointed Gromov-Hausdorff topology.
By definition all $D'_n$ are unbounded, so $D'$ is also unbounded.

Let $\theta_p$ and $\theta_{p'}$ be the cone angles at $p$ and $p'$,
respectively, in $M$, so that $\theta_p,\theta_{p'}\leq \theta$.
We now consider $D'$ as a convex polyhedron in Euclidean space $\R^3$,
with two edges $e$ and $e'$ corresponding to $p$ and $p'$, respectively.
Let $H_1, H_2$ be the two half-planes bounded by $p$ at angle
$\theta_p/2$ with $s$, and let $H'_1, H'_2$ be the two half-planes 
bounded by $p'$ at angle $\theta_{p'}/2$ with $p'$. 
Then $H_1, H_2, H'_1, H'_2$ are faces of $D'$, so that $D'$ is contained
in $D''$, the intersection of the half-spaces bounded by the four planes 
containing $H_1,H_2,H'_1, H'_2$ and containing $s$.

Suppose that $e$ and $e'$ are not parallel. Then $D''$ has at most one end, so 
that $D'$ has also at most one end. This is clearly impossible since all $M_n$ are
quasifuchsian manifolds with particles, so that all $D_n$ have two ends.
Therefore, $e$ and $e'$ are parallel. For the same reason any other edge of $D'$ which corresponds or not to a cone singularity of $M$ has to be
parallel to $e$ and to $e'$. So $D'$ is invariant under translations 
parallel to $e$, that it, it is the product by $\R$ of a polygon
$\pi$ in a plane orthogonal to $e$. It follows that $M'$ is also invariant
by translation parallel to $p$.

Consider the regular part $M'_{reg}$ of $M'$. 
Since $M'_{reg}$ is a Euclidean manifold, its holonomy representation $Hol$ is a 
morphism from $\pi_1(M'_{reg})$ to $\isom(\R^3)=\R^3\rtimes O(3)$. Since $M'$
is invariant under translations parallel to $p$, $Hol$ actually takes values in 
$\R^3\rtimes O(2)$. We consider the morphism $Hol':\pi_1(M'_{reg})\rightarrow 
\R^2\rtimes O(2)=\isom(\R^2)$ obtained by projecting the translation component
of each element on the plane orthogonal to $p$. 

Then $Hol'$ is the holonomy representation of a 
3-dimensional %% ajoute
Euclidean cone-manifold $M'_{par}$.
%which is the quotient of $M_{reg}$ under translation in the invariant direction.
(We do not discuss whether $Hol=Hol'$ in all cases.)
By construction $M'_{par}$ contains a complete surface orthogonal to the singular
locus, say $S$. This surface has at least two singular points, and each of its
singular points has angle equal to the angle of the
corresponding cone singularity of $M'$, so those angles are less than $\theta$.
Since $\theta<\pi$ it follows from the Gauss-Bonnet formula that $S$ is 
homeomorphic to a sphere, and that it has at most three cone singularities. 

Still by construction, the fundamental group of $M'_{reg}$ surjects to the fundamental
group of $M'_{par}$. However this surjection is actually an isomorphism, since
otherwise an element of $\pi_1(M'_{reg})$ would act trivially on $S$, which means
that it would act on $M'_{reg}$ by translations parallel to the invariant direction,
and this is impossible since $M'_{reg}$ is non-compact. 

This shows that $M'$ is homeomorphic to the product of a sphere by a line, and that it has three
cone singularities. Therefore this is also true of all $M_n$ for $n$ large
enough. But this is impossible since the $M_n$ are quasifuchsian manifolds
with particles, and the definition explicitly excludes manifolds homeomorphic
to the product of a sphere by a line with three singularities.
\end{proof}

We now call $\epsilon_0>0$ the number $\epsilon$ associated by the previous two
lemmas to the maximum of the $\theta_i$, and $\rho_0$ the corresponding
value of $\rho$.  

%% file: conebend3.tex
\newcommand{\eps}{\varepsilon}
\newcommand{\Hp}{\mathbb{H}}

\section{Compactness statements}
\label{se:compact-lamination}

%% !TEX root =conebend.tex

\subsection{Main statement.}    \label{main}

The main goal of this section is to prove the following compactness lemma. 

\begin{lemma} \label{lm:compact}
Let $M_n$ be a sequence of quasifuchsian manifolds with particles with the same topological type and converging angles.
Let $\lambda_n$ be the measured bending 
laminations on the boundary of the convex core of $M_n$. 
Suppose that $\lambda_n\rightarrow \lambda_\infty$, where $\lambda_\infty$ satisfies 
the hypothesis of Theorem \ref{tm:general}. 
Then, after taking a subsequence, $M_n$ converges to a quasifuchsian
 manifold with particles with the common topological type, the limit particles and  measured bending lamination $\lambda_\infty$.
\end{lemma}

Let us explain the definitions used in this statement.
The {\em topological type} of a quasifuchsian manifold with particles $M$ has the form $(S,x_1,\cdots, x_{n_0})$ where $S$ is a compact surface with genus at least $2$ and $x_1,\cdots, x_{n_0}$ are distinct points on $S$. A quasifuchsian manifold with particles $M$ has topological type $(S,x_1,\cdots, x_{n_0})$ if $M$ is isometric to the product $S\times\R$ endowed with a complete hyperbolic metric with cone singularities of angles $\theta^i\in (0,\pi)$ on the lines $\{x_i\}\times\R$.

Consider a sequence of quasifuchsian manifolds $M_n$ with particles with the same topological type $(S,x_1,\cdots, x_{n_0})$. Denote by $\theta^i_n\in (0,\pi)$ the cone angles of the metric of $M_n$ on $\{x_i\}\times\R$. Then the sequence $M_n$ has {\em converging angles} if and only if $\theta^i_n$ converges in $(0,\pi)$ when $n$ goes to $\infty$ for any $i\leq n_0$.

Notice that since $\lambda_n$ converge to $\lambda_\infty$, $\lambda_n$ is eventually non trivial. 
In particular the manifolds $M_n$ are not fuchsian (i.e. their convex core is not a surface) except maybe for finitely many of them. Throughout %% of 
this section, when we consider a quasifuchsian manifold with particle, we will assume that it is not fuchsian so that its convex core is a $3$-dimensional manifold with boundary.

%% paragraphe ajoute
The convex core $C(M)$ of $M$ is homeomorphic to $S\times I$. Thus $\partial C(M)$ is homeomorphic to $S\sqcup S$ and each copy of $S$ in this union has $k$ marked points $x_1,\cdots, x_{n_0}$ corresponding to the endpoints of the particles. The measured bending 
%% laminations 
lamination
on the boundary of $C(M)$ is an element of the space $\cML_{S,n_0}\times \cML_{S,n_0}$ of measured laminations on two copies of $S$ with $n_0$ marked points. The space $\cML_{S,n_0}$ is endowed with the topology of weak-$*$ convergence of measures on
compact transversals and $\cML_{S,n_0}\times \cML_{S,n_0}$ is endowed with the product topology. In simple terms, we can fix a finite but sufficiently large set of curves $c_i$ which are either closed curves or segments between two singular points, then two measured laminations are close if and only if their intersection with each of the $c_i$ are close.  
%% fin paragraphe ajoute

Although Lemma \ref{lm:compact} is a generalisation of the "Lemme de fermeture" of \cite{bonahon-otal}, 
the proof is very different. The reason is that the two main ingredients of the proof in 
Bonahon-Otal's paper are Culler-Morgan-Shalen compactification of the character variety 
by actions on $\R$-trees and the covering Theorem of Canary. Since both these results 
hardly extend to manifolds with particles we had to use different arguments. 
Since our proof also works without particles, we get a new proof of the main result of 
\cite{bonahon-otal}.

\subsection{A finite cover argument} \label{ssc:finite}

We work under the assumption that the cone angle around each singularity is less than $\pi$. 
This assumption guarantees that the singularities are never too close to each other, 
see Lemma \ref{lm:dist3},
and that the boundary of the convex core is well defined and is orthogonal to the singularities. 
On the other hand, cone singularities with cone angles less than $\pi$ can be viewed has singularities 
with concentrated positive curvature. 
But some of the results we will use are easier to prove when the curvature is negative. 
To overcome this difficulty, we will use a branched cover for which the cone angles are all greater than $2\pi$.

Let $M$ be a quasifuchsian manifold with particles. A branched cover $\overline{M}\rightarrow M$ branched 
along the singularities is {\em negatively curved} if the cone angles around the singularities of the metric induced on $\overline{M}$ %% or 
are all greater than $2\pi$. We call $\overline{M}$ {\em a negatively curved branched cover of $M$}.

This name comes from the fact that a singularity with cone angle greater than $2\pi$ can be viewed as a set of concentrated negative curvature. 
More precisely $\overline{M}$ can be approximated by Riemannian manifolds with curvature bounded above by $-1$ (in the bilipschitz topology). 
It follows that $\overline{M}$ has properties of negatively curved manifolds, in particular the uniqueness of the geodesic segment joining two given points in a given homotopy class.

We will construct such branched covers for sequences.
Consider a sequence of quasifuchsian manifold with particles $M_n$ with the same topological type $(S,x_1,\cdots, x_{n_0})$ (as defined in the preceding section). We denote by $g_n$ the metric of $M_n$ and by $\theta^i_n$ the cone angles of $g_n$ on $\kappa_i$. Assume that the sequence $M_n$ has converging angles, namely $\theta^i_n$ converge to some $\theta^i\in (0,\pi)$ for any $i\leq n_0$.

For each singularity $x_i$, we choose an integer $k_i$ such that $\frac{2\pi}{k_i}$ 
is less than the angle $\theta^i$ (the limit of $\theta^i_n$). 
The surface $S$ with cone angle $\frac{2\pi}{k_i}$ at the point $x_i$ is a hyperbolic orbifold. 
As such it has a manifold cover $h:\bar S\rightarrow S$ which is a branched cover 
so that the lifts of the point $x_i$ have a branching index equal to $k_i$. 
The branched cover $h:\bar S\rightarrow S$ extends naturally to a branched cover 
$h:\bar S\times I\rightarrow S\times I$.

For a fixed $n$, we have the metric $g_n$  on $S\times\R$ with cone singularities $\theta^i_n$ along $\{x_i\}\cap\R$. If we pull back $g_n$ using the map $h$, we get a hyperbolic metric $\bar g_n$ with cone singularities on $\bar S\times \R$ for which the covering transformations are isometries. 
Let $\overline{M}_n=(\bar S\times\R,\bar g_n)$ be the manifold with cone singularities thus obtained. By the choice of $\{k_1,...,k_{n_0}\}$, for $n$ large enough, we have $k_i\theta^i_n\geq 2\pi$, hence the cone angle of $\bar g_n$ around each singularity of $\overline{M}$ is at least $2\pi$. Thus for $n$ large enough, $\overline{M}_n$ is a negatively curved branched cover of $M_n$ and the topological type of $\overline{M}_n$ does not depend on $n$.

\subsection{Pleated annuli}	\label{pa}

A technical device which will be useful later on is a simplicial annulus bounded by two given curves. As was mentioned above, when we consider a quasifuchsian manifold with particle, we assume that it is not fuchsian.

Let us first fix some notations. 
We %%have 
consider a quasifuchsian %% ajoute "a"
manifold with particles $M$ with topological type $(S,x_1,\cdots, x_{n_0})$. 
We denote by $g$ the complete hyperbolic metric  with cone singularities of $M$ and by $ C(M)$ the convex core of $M$. 
We will use a negatively curved branched cover $\overline{M}$ of $M$ (as defined in the previous section). 
The construction of such a cover is explained above for a sequence $M_n$, here we take the constant sequence, $M_n=M$ for any $n$, to define $\overline{M}$. We denote by $ C(\overline{M})\subset\overline{M}$ the preimage of the convex core $ C(M)$ of $M$ under the covering projection, by $\bar\lambda$ its bending measured geodesic lamination and by $\bar m$ the induced metric on $\partial  C(\overline{M})$. We will use these notations throughout this section.

Now let us construct our simplicial annulus.

\begin{lemma} \label{lm:annulus}
Let $M$ be a quasifuchsian manifold with particles and $\overline{M}$ a negatively curved branched cover of $M$. 
Let $\bar d,\bar d'$ be homotopic simple closed geodesics, 
respectively on the upper and on the lower boundary component of $ C(\overline{M})$. 
There exists an immersed annulus $\bar A$ in %% $ C(\bar  N)$ 
$C(\bar  M)$ 
bounded by %% $\bar d \cup \bar d'\subset\partial C(N)$ 
$\bar d \cup \bar d'\subset\partial C(\bar M)$ 
such that the metric induced on $\bar A$ by $g_n$ is a hyperbolic metric 
with cone singularities with angles at least $2\pi$. 
The area of $\bar A$ is at most $\max\{l_{\bar m} (\bar d)+l_{\bar m} (\bar d'), i(\bar\lambda,\bar d)+i(\bar\lambda,\bar d')\}$.
\end{lemma}

\begin{proof}
Let us specify that $l_{\bar m} (\bar d)$, resp. $l_{\bar m} (\bar d')$, is the length of $\bar d$, resp. $\bar d'$, with respect to the metric $\bar m$ induced by $\bar g$ on $\partial  C(\overline{M})$.
 
Since $\bar d$ and $\bar d'$ are disjoint homotopic simple closed curves, there is an embedded annulus $\bar A\subset  C(\overline{M})$ with $\partial\bar A=\bar d\cup\bar d'$. If the bending laminination of $ C(\overline{M})$ intersects $\bar d$ and $\bar d'$ finitely many times, then $\bar d$ and $\bar d'$ are piecewise geodesics. If not, we approximate them by piecewise geodesic curves and work on the approximates. Consider a triangulation $T$ of $\bar A$ whose vertices are all contained in $\bar d\cup\bar d'$ and such that any vertex of $\bar d$ and $\bar d'$ (when considered as piecewise geodesics) is a vertex of $T$. As we have said before, in $\overline{M}$, there is a unique geodesic segment joining $2$ given points in a given homotopy class. It follows that we can change $\bar A$ by a homotopy so that  each edge of $T$ is a geodesic segment in $ C(\overline{M})$. 
Next, for each triangle $T_i$ of $T$, we choose a vertex $v$ and we substitute $T_i$ by the geodesic cone from $v$ to the edge $e_v$ of $T_i$ not containing $v$. 
This geodesic cone is the union of the geodesic segments joining $v$ to the edge $e_v$ of $T_i$ (the homotopy class of such segment is defined by the corresponding segment of $T_i$). Again the existence of this cone follows from the uniqueness of geodesic paths. From now on we denote this cone by $T_i$. By construction, it is a locally ruled surface and as such has negative curvature:

\begin{claim}		\label{triangle}
Let $M$ be a hyperbolic manifold with cone singularities with cone angles bigger than $2\pi$. Given a point $v\in M$ and a geodesic segment $e_v\subset M$, a geodesic cone $T_i\subset M$ from $v$ to $e_v$ is an union of polygons with curvature $-1$.  Furthermore, the sum of the angles of the polygons meeting at an interior vertex is at least $2\pi$.
\end{claim}

\begin{proof}
The surface $T_i$ meets the singular locus $\overline{M}_s$ of $\overline{M}$ along segments and at points. For each component $\kappa$ of $\overline{M}_s\cap T_i$ we consider the two extremal segments joining $v$ to $e_v$ and intersecting $\kappa$. Doing this for each component of $\overline{M}_s\cap T_i$, we get a family of segments which are geodesic for the metric of $(\overline{M},\bar g)$ and hence for the induced metric on $T_i$. We add the components of $\overline{M}_s\cap T_i$ which are segments to this family and get a new family of geodesic segments. The closure of each complementary region is a polygon, i.e. a disc with piecewise geodesic boundary (see Figure \ref{figannulus}). By construction each such polygon is a locally ruled surface in $\Hp^3$ hence it has curvature $-1$. Thus we have proved the first sentence of this Claim.

\begin{figure}[hbtp]            \label{figannulus}
\psfrag{a}{$\overline{M}_s\cap T_i$}
\centerline{\includegraphics{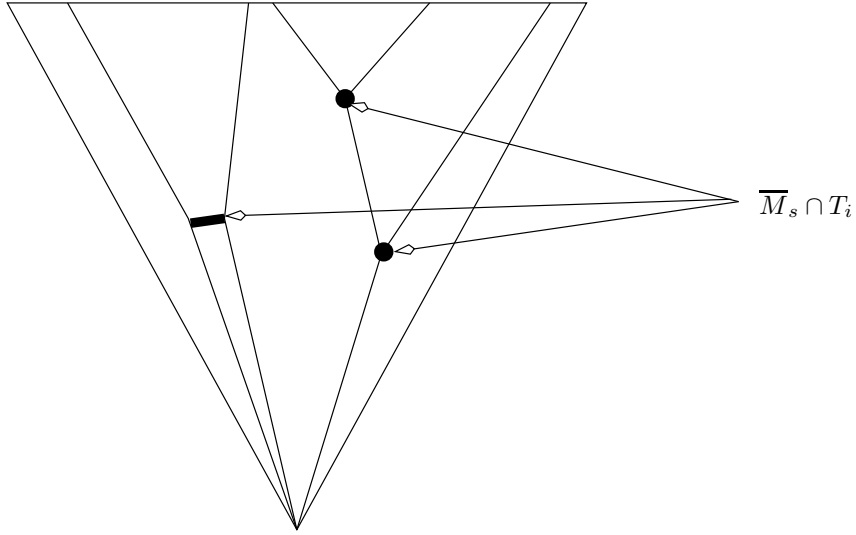}}
\caption{Decomposition of $T_i$ into hyperbolic polygons} 
\end{figure}

%Since no two components of the singular locus intersect each other,  any vertex of any polygons lies in at most one component of the singular locus. Furthermore,
By construction, given an interior vertex $v$ of this decomposition into polygons, there is a geodesic segment (for the metric of $(\overline{M},\bar g)$) which passes through $v$. On each side of this segment, the sum of the angles of the polygons has to be at least $\pi$. Thus, we can conclude that the sum of the angles of the polygons around $v$ is at least $2\pi$.
\end{proof}

\indent
We change the annulus $\bar A$ so that it is a union of geodesic cones as described in Claim \ref{triangle}. Thus %% that 
the induced metric is hyperbolic with cone singularities with angles greater than $2\pi$. By the %% ajoute "the"
Gauss-Bonnet Formula, the area of $\bar A$ is at most the bending of $\partial\bar A$, namely ${\rm Area}(\bar A)\leq i(\partial\bar A,\bar\lambda_n)=i(\bar d,\bar\lambda_n)+i(\bar d',\bar\lambda_n)$.

It remains to prove that ${\rm Area}(\bar A)\leq l_{\bar m} (\bar d)+l_{\bar m} (\bar d')$. By construction, $\bar A$ is a union of triangles $T_i$ such that one edge of each $T_i$ lies in $\bar d\cup\bar d'$ and by Claim \ref{triangle}, the induced metric on each such triangle is a hyperbolic metric with cone singularities with angles greater than $2\pi$. It follows that the induced metric can be approximated by Riemannian metrics with curvature at most $-1$. 
Let $T_i^h$ be a hyperbolic triangle (i.e. a geodesic triangle in $\Hp^2$) such that the length of the edges of $T_i^h$ are the same as the length of the edges of $T_i$. Since the induced metric on $T_i$ as curvature at most $-1$, we have ${\rm Area}(T_i)\leq {\rm Area}(T_i^h)$. On the other hand the area of a hyperbolic triangle is less than the length of any of its edges (see \cite[Lemma 9.3.2]{thurston-notes}). It follows that ${\rm Area}(T_i)$ is less than the length of any of its edges, in particular it is less than the length of the edge of $T_i$ lying in $\bar d\cup \bar d'$. Since this holds for all the triangles composing $\bar A$, we have ${\rm Area}(\bar A)\leq l_{\bar m} (\bar d)+l_{\bar m} (\bar d')$.
\end{proof}

\subsection{Long geodesics in $M$.} \label{ssc:34}

In this section, we will show that, under the hypothesis of Lemma \ref{lm:compact}, 
the induced metrics on $\partial C(M_n)$ are bounded. In order to do that we will show that 
if some geodesic is long in the boundary of $ C(M)$, then the boundary of some annulus is almost not bent or the bending lamination tends to have a leaf with a weight greater than or equal to $\pi$.  Since this would contradict the conditions on $\lambda_\infty$, it will follow that any given simple closed curve on $\partial C(M_n)$ has bounded length. As earlier, when we consider a quasifuchsian manifold with particle, we assume that it is not fuchsian.\\

\indent
Throughout this section we use the following notations. %%:
We have a sequence of quasifuchsian manifolds with particles $M_n$ with the same topological type $(S,x_1,\cdots, x_{n_0})$. 
We denote by $(g_n)_{n\in \N}$ the metric of $M_n$. We assume that the sequence $M_n$ has converging angles, 
namely $\theta^i_n$ converge in $(0,\pi)$ for any $i\leq n_0$. 
Since $M_n$ is quasifuchsian (with particles), $C(M_n)$ is homeomorphic to $S\times I$. 
We denote by $S$ and $S'$ the two components of $\partial C(M_n)$, $S'$ is homeomorphic to $S$. 
Let $m_n$ be the metric defined on $S\sqcup S'$ by the identification with $\partial  C(M_n)$ endowed with the metric induce by the $g_n$-length of paths. This metric $m_n$ is a hyperbolic metric with cone singularities of angles $\theta^i_n$ at the points $\{x_i\}\in S$ and $\{x'_i\}\in S'$. We denote by $\lambda_n\in \cML_{S,n_0}\times \cML_{S,n_0}$ the bending measured geodesic lamination of $\partial C(M_n)$.

We will make use of the branched covers $\overline{M}_n$ defined in Section \ref{ssc:finite}. Recall that, for $n$ large enough, $\overline{M}_n$ is a negatively curved branched cover of $M_n$ and that the topological type of $\overline{M}_n$ does not depend on $n$. We denote by $ C(\overline{M}_n)$, $\bar m_n$ and $\bar \lambda_n$ the preimages of $ C(M_n)$, $m_n$ and $\lambda_n$ under the covering map $\overline{M_n}\rightarrow M_n$.\\

\indent
We will use the next lemma to prove that, under the right hypothesis on $\lambda_\infty$, the induced metric on $\partial C(M_n)$ is bounded.

\begin{lemma}   \label{annulusorpi}
Let $M_n$ be a sequence of quasifuchsian manifolds with particles with the same topological type and converging angles. Assume that $(\lambda_n)$ converges to $\lambda_\infty$ (without any hypothesis on $\lambda_\infty$). 
Consider a simple closed curve $d\subset S$. Let $d_n\subset \partial C(M_n)$ be the closed 
$m_n$-geodesic freely homotopic to $d$. 
If $l_{m_n}(d_n)\longrightarrow\infty$, then either $\lambda_\infty$ contains a leaf with a weight greater than or equal to $\pi$, 
or there is a sequence of essential annuli $E_n$ such that $i(\lambda_n,\partial E_n)\longrightarrow 0$.
\end{lemma}

\begin{proof}
Let $S'$ be the other boundary component of $ C(M)$ (i.e. not $S$), 
and let $d'_n$ be the closed $m_n$-geodesic freely homotopic to $d$ lying in $S'$. Let $\overline{M}_n$ be a negatively curved branched cover of $M_n$ so that for $n$ large enough the topological type of $\overline{M}_n$ does not depend on $n$.
 Let $\bar d_n$ and $\bar d'_n\subset\partial  C(\overline{M}_n)$ be homotopic lifts of $d_n$ and $d'_n$ respectively under the covering projection $\overline{M}\rightarrow M$. The preimage $\bar\lambda_n\in{\cal ML}(\partial\overline{M})$ of $\lambda_n$ is the bending measured lamination of $ C(\overline{M}_n)$. Furthermore $\bar\lambda_n$ converges %% converge 
to the preimage $\bar\lambda_\infty$ of $\lambda_\infty$.\\
\indent
First we will show that if $\bar d_n$ is long compared to the area of an annulus $\overline{A}_n$ bounded by $\bar d_n\cup\bar d'_n$ then there are shortcuts in $\overline{A}_n$. Namely $\bar d_n\cup\bar d'_n$ contains points which are close to each other in $\overline{M}_n$ but far
in $\bar d_n\cup\bar d'_n$. This can happen for instance if $d_n$ and $d'_n$ are close to each other in $ C(M_n)$.

\begin{claim}           \label{petitarc}
Let $M_n$ be a sequence of quasifuchsian manifolds with particles and let $\overline{M}_n$ be a negatively curved branched cover of $M_n$ such that the topological type of $\overline{M}_n$ does not depend on $n$. Assume that %%$(\bar\lambda_n)$ 
$(\bar\lambda_n)_{n\in \N}$ 
converges to $\bar\lambda_\infty$. Let $\bar d\subset\bar S$ be a simple closed curve and denote by $\bar d_n\subset\bar S\subset \partial  C(\overline{M}_n)$ the simple closed $\bar m_n$-geodesic in the homotopy class of $\bar d$. 

If $l_{\bar m_n}(d_n)\longrightarrow\infty$, then there is a $\bar m_n$-geodesic arc $\bar k_n\subset  C(\overline{M}_n)$ 
such that $\ell_{\bar m_n}(\bar k_n)\longrightarrow 0$ and that either $\bar k_n$ joins the two components of $\partial C(\overline{M}_n)$ 
or the $\bar m_n$-geodesic arc %% $\bar\kappa_n\in \bar d_n$ 
$\bar\kappa_n\subset \bar d_n$ 
in the homotopy class of $\bar k_n$ relative to its boundary satisfies $\ell_{\bar m_n}(\bar \kappa_n)\longrightarrow\infty$.
\end{claim}

\begin{proof}
Denote by $\bar d'_n$ the closed $\bar m_n$-geodesic lying in $\bar S'\subset \partial  C(\overline{M}_n)$ that is homotopic to $\bar d_n$ in $ C(\overline{M}_n)$.
Consider the annulus $\bar A_n$ with $\partial\bar A_n=\bar d_n\cup\bar d'_n$ 
that was constructed in Lemma \ref{lm:annulus}. 
Since $l_{\bar m_n}(\bar d_n)\longrightarrow\infty$, there is $\eps_n\longrightarrow 0$ and a segment $\bar s_n\subset\bar d_n$ such that 
$l_{\bar m_n}(\bar s_n)\longrightarrow\infty$ and $i(\bar s_n, \bar\lambda_n)\leq\eps_n$. 
Let $\bar t_n\subset  C(\overline{M}_n)$ be the $\bar g_n$-geodesic segment homotopic to $\bar s_n$ relative to its endpoints. 
Since $\bar s_n$ is almost not bent, its length is very close to 
the length of $\bar t_n$ (see \cite[Lemme A2]{lecuire}). 
In particular, $l_{\bar g_n}(\bar t_n)\longrightarrow\infty$. 
Furthermore, for the same reason, any point in $\bar s_n$ is close to $\bar t_n$. 
Namely there is $\eta_n=\eta(\varepsilon_n)\longrightarrow 0$ such that 
for any point $\bar z_n\subset\bar s_n$, 
there is $\bar x_n\subset\bar t_n$ with $d_{\bar g_n}(\bar x_n,\bar z_n)\leq\eta_n$ 
(see \cite[Affirmation A3]{lecuire}).

Since $(\bar\lambda_n)_{n\in \N}$ converges to $\bar\lambda_\infty$, then the bending $i(\bar\lambda_n,\partial\bar A_n)$ of $\partial\bar A_n$ converges. By Lemma \ref{lm:annulus} the area of $\bar A_n$ is bounded. Now, in $\bar d_n$, we replace $\bar s_n$ by $\bar t_n$. By the previous paragraph, we can still consider the annulus $\bar A_n$ and its area is bounded. For any point in $\bar t_n$ that is at distance at least $\frac{1}{3}l_{\bar m_n}(\bar t_n)$ from $\partial \bar t_n$, we consider in $\bar A_n$ an arc orthogonal to $\bar t_n$ that either hits $\partial\bar A_n$ 
at distance less than $\eta_n$ from its basepoint %% ajout'e
or has length $\eta_n$ ($\eta_n$ will be specified later on). Let $\bar Z_n\subset \bar A_n$ be the union of those arcs that have length $\eta_n$ and let $\bar z_n$ be the union of their starting points (i.e. their intersection with $\bar t_n$). The set $\bar Z_n$ is embedded and its area is the same as the area of a strip of length $\ell_{\bar m_n}(\bar z_n)$ and width $\eta_n$.  Notice that since the singularities of $\bar A_n$ have cone angles at least $2\pi$, the area of this strip at least the area of a hyperbolic strip with the same length and width,  i.e. it is at least $\ell_{\bar m_n}(\bar z_n)\sinh (\eta_n)$. Let $K$ be a number larger than %% an upper bound for 
the area of $\bar A_n$. 
Taking $\eta_n$ such that $\sinh(\eta_n) >\frac{3K}{\ell_{\bar m_n}(\bar t_n)}$, we get  
$$ K\geq {\rm Area}(\bar Z_n)\geq  \ell_{\bar m_n}(\bar z_n)\sinh (\eta_n)>K\frac{3\ell_{\bar m_n}(\bar z_n)}{\ell_{\bar m_n}(\bar t_n)}~. $$
Hence $\ell_{\bar m_n}(\bar z_n)<\frac{1}{3} \ell_{\bar m_n}(\bar t_n)$. It follows that there exists an arc with length less than $\eta_n$ orthogonal to $\bar t_n$ whose starting point $\bar x'_n\subset\bar t_n$ is at a distance at least $\frac{1}{3}l_{\bar m_n}(\bar t_n)$ from $\partial \bar t_n$ (distance measured on $\bar t_n$) and which hits $\partial\bar A_n$ on a point $\bar y_n\subset\partial\bar A_n$.

As we have seen in the previous paragraph there is a point $\bar x_n\subset\bar k_n$ very close to $\bar x'_n$. It follows that $\bar x_n\subset\bar d_n$ and $\bar y_n\subset \bar d_n\cup\bar d'_n$ are joined in $ C(\overline{M}_n)$ by an arc $\bar k_n$ satisfying $\ell_{g_n}(\bar k_n)\longrightarrow 0$.

If $\bar y_n\in \bar d'_n$ then we are done. Otherwise $\bar x_n$ and $\bar y_n$ both lie in $\bar d_n$. By construction $\bar k_n$ lies in an annulus connecting $\bar d_n$ to $\bar d'_n$. It follows that there is a $\bar m_n$-geodesic arc $\bar\kappa_n\subset \bar d_n$ that is homotopic to $\bar k_n$ relative to $\{\bar x_n\}\cup\{\bar y_n\}$. Since $\bar x_n$ is at distance at least $\frac{l_{\bar m_n}(\bar t_n)}{3}$ from the points in $\partial\bar s_n$, we have 
$$ \ell_{\bar m_n}(\bar \kappa_n)\geq \frac{l_{\bar m_n}(\bar t_n)}{3}\longrightarrow\infty~. $$
\end{proof}

Consider the points $\bar x_n$ and $\bar y_n$ constructed in Claim \ref{petitarc} and extract a subsequence such that either $\bar y_n\in\bar d_n$ for any $n$ or $\bar y_n\in\bar d'_n$ for any $n$. 
We will show below that if $\bar y_n$ lies in $\bar d_n$ then 
$\lambda_\infty$ has a leaf with a weight greater than or equal to $\pi$ and that if $\bar y_n$ lies in $\bar d'_n$ then there is a sequence of essential annuli $E_n\subset M$ such that $i(\lambda_n,\partial E_n)\longrightarrow 0$.

In the next step we are going to construct $\bar m_n$-geodesic loops based at $\bar x_n$ and $\bar y_n$ that are almost not bent.

\begin{claim}   \label{anneauplat}
%%Consider a $\bar m_n\subset\partial C(\overline{M}_n)$-geodesic loop $\bar l_n$ based at $\bar x_n$.
Let $M_n$ be a sequence of quasifuchsian manifolds with particles and let $\overline{M}_n$ be a negatively curved branched cover of $M_n$ such that the topological type of $\overline{M}_n$ does not depend on $n$. Consider $2$ points $\bar x_n,\bar y_n\in\partial C(\overline{M}_n)$ away from the singularities and a $\bar g_n$-geodesic arc $\bar k_n\subset C(\overline{M}_n)$ joining $\bar x_n$ to $\bar y_n$ such that $\ell_{\bar g_n}(\bar k_n)\longrightarrow 0$ and that either $\bar x_n$ and $\bar y_n$ lie on different components of $\partial C(M_n)$ or there is a $\bar m_n$-geodesic arc $\bar\kappa_n\in \partial C(\overline{M}_n)$ that is homotopic to $\bar k_n$ relative to its boundary $\{\bar x_n\}\cup\{\bar y_n\}$ and that satisfies $\ell_{\bar m_n}(\bar \kappa_n)\longrightarrow\infty$.
 
Consider a loop $\bar l_n$ on %% $ C(\overline{M}_n)$ 
$ \partial C(\overline{M}_n)$ 
based at $\bar x_n$, which is geodesic for $\bar m_n$ 
(except at $\bar x_n$). 
Let $\bar f_n\subset \partial C(\overline{M}_n)$ be the $\bar m_n$-geodesic loop based at $\bar y_n$ 
that is homotopic to $\bar l_n$. 
%%Assume that $\ell_{\tilde f_n}(\bar g_n)$ is bounded. 
Assume that $\ell_{\bar m_n}(\bar l_n)$ is bounded. 
Then the bending of $\bar l_n$ and $\bar f_n$ tends to $0$, namely $i(\bar l_n, \bar \lambda_n)\longrightarrow 0$ and $i(\bar f_n, \bar \lambda_n)\longrightarrow 0$.
\end{claim}

\begin{proof}
When saying that $\bar x_n$ and $\bar y_n$ are away from the singularities we mean that there is a uniform upper bound on %% fro
their distance to the singular locus of $\overline{M}_n$.

Let $\widetilde{M}_n$ be the universal cover of $\overline{M}_n$, it is a simply connected hyperbolic $3$-manifold with cone singularities.  Let $C(\widetilde{M}_n)$ be the lift of $ C(\overline{M}_n)$ to $\widetilde{M}_n$. Let $\tilde l_n$, $\tilde k_n$, $\tilde x_n$, $\tilde y_n$ be lifts of $\bar l_n$, $\bar k_n$, $\bar x_n$ and $\bar y_n$ with $\tilde x_n\in \tilde l_n$ and $\tilde x_n\cup\tilde y_n=\partial\tilde k_n$. 
The point $\bar x'_n=\partial\tilde l_n\setminus \tilde x_n$ is the image of $\bar x_n$ under a covering transformation.
Consider the $\tilde m_n$-geodesic arc $\tilde f_n\subset\partial C(\widetilde{M}_n)$ 
joining $\bar y_n$ to its image $\bar y'_n$ under this covering transformation.

Let us first assume that there are no singularities in $\overline{M}_n$. Then $\widetilde{M}_n$ is isometric to $\Hp^3$ and we choose the isometry so that $\tilde x_n$ is identified with a fixed point of $\Hp^3$ (independantly of $n$). Let $\Pi(\tilde x_n)$ be a support plane for $C(\widetilde{M}_n)$ at $\tilde x_n$, namely a totally geodesic plane that intersects $C(\widetilde{M}_n)$ only along $\partial C(\widetilde{M}_n)$ and contains $\tilde x_n$. Up to moving $\bar x_n$ slightly, we may assume that it is disjoint from $\bar\lambda_n$ so that there is only one support plane at $\tilde x_n$. The convex set $C(\widetilde{M}_n)$ lies in a half-space $E(\tilde x_n)$ bounded by $\Pi(\tilde x_n)$. Similarly let $\Pi(\tilde y_n)$ be a support plane at $\tilde y_n$ and let $E(\tilde y_n)$ be the half-space bounded by $\Pi(\tilde y_n)$ that contains $C(\widetilde{M}_n)$.

If $\bar x_n$ and $\bar y_n$ lie on different components of $\partial C(M_n)$ then $\Pi(\tilde x_n)$ and $\Pi(\tilde y_n)$ are disjoint. Otherwise, since $\ell_{\bar m_n}(\bar \kappa_n)\longrightarrow\infty$, either $\Pi(\tilde x_n)$ and $\Pi(\tilde y_n)$ are disjoint or their intersection goes to $\infty$ with $n$ (namely the sequence $\Pi(\tilde x_n)\cap\Pi(\tilde y_n)$ lies outside larger and larger compact sets in $\Hp^3$). Since $\ell_{\bar g_n}(\bar k_n)\longrightarrow 0$, $\tilde y_n$ converges to $\tilde x_n$ (viewed as a fixed point in $\Hp^3$) and, up to extracting a subsequence, $\Pi(\tilde x_n)$ and $\Pi(\tilde y_n)$ converges to the same plane $\Pi_\infty$ in $\Hp^3$. Furthermore $E(\tilde x_n)$ converges to a half-space $E(\tilde x_\infty)$ bounded by $\Pi_\infty$ and $E(\tilde y_n)$ converges to the other half-space $E(\tilde y_\infty)$ bounded by $\Pi_\infty$. 

Let $\Pi(\tilde x'_n)$ be a support plane at $\tilde x'_n$. Since $\ell_{\bar m_n}(\bar f_n)$ is bounded, up to extracting a subsequence, $\tilde x'_n$ converges in $\Hp^3$. Again since  either $\bar x'_n$ and $\bar y_n$ lie on different components of $\partial C(M_n)$ or $d_{\tilde m_n}(\tilde x'_n,\tilde y_n)\longrightarrow\infty$, $\Pi(\tilde y_n)\cap\Pi(\tilde x'_n)$ either is empty or goes to infinity. It follows that $\Pi(\tilde x'_n)$ also converges to $\Pi_\infty$ and that $E(\tilde x'_n)$ converges to $E(\tilde x_\infty)$. The external dihedral angle between $\Pi(\tilde x_n)$ and $\Pi(\tilde x'_n)$ is an upper bound for  $i(\bar l_n, \bar \lambda_n)$, hence  $i(\bar l_n, \bar \lambda_n)\longrightarrow 0$.

It remains to show that %% ajoute "that"
$i(\bar f_n, \bar \lambda_n)\longrightarrow 0$.

By construction, $d_{\tilde g_n}(\tilde x'_n,\tilde y'_n)=d_{\tilde g_n}(\tilde x'_n,\tilde y'_n)=\ell_{\bar g_n}(\bar k_n)\longrightarrow 0$, hence $\Pi(\tilde x'_n)$ and $\Pi(\tilde y'_n)$ converge to the same plane $\Pi_\infty$ in $\Hp^3$. If $\bar x_n$ and $\bar y_n$ lie on different components of $\partial C(M_n)$ then $\Pi(\tilde x'_n)$ and $\Pi(\tilde y'_n)$ are disjoint. It follows that $E(\tilde y'_n)$ converges to $E(\tilde y_\infty)$ while $E(\tilde x'_n)$ converges to $E(\tilde x_\infty)$ which implies $i(\bar f_n, \bar \lambda_n)\longrightarrow 0$ as above. If $\bar x_n$ and $\bar y_n$ lie on the same component of $\partial C(M_n)$ then $d_{\tilde m_n}(\tilde x'_n,\tilde y'_n)=d_{\tilde m_n}(\tilde x_n,\tilde y_n)\longrightarrow\infty$. It follows that $E(\tilde y'_n)$ converge to $E(\tilde y_\infty)$ which again implies $i(\bar f_n, \bar \lambda_n)\longrightarrow 0$.

When $\overline{M}_n$ have singularities, we cannot define support planes, but we can define local support planes at points which are disjoint from the singularities. Thus we can locally use the same arguments as in the non singular case, leading to the same conclusion.
\end{proof}

Let us choose for $\bar l_n$ a shortest $\bar m_n$-geodesic loop based at $\bar x_n$. 
Since the area of $(\bar S\sqcup\bar S',\bar m_n)$ is bounded, 
there is a constant $Q>0$ such that $\ell_{\bar m_n}(\bar l_n)\leq Q$. 
By Claim \ref{anneauplat}, we have $i(\bar l_n, \bar \lambda_n)\longrightarrow 0$. 
Let $\bar f_n$ be the $\bar m_n$-geodesic loop based at $\bar y_n$ that is homotopic to $\bar l_n$.
By Claim \ref{anneauplat}, we have $i(\bar f_n, \bar \lambda_n)\longrightarrow 0$. 
Since $\bar l_n$ and $\bar f_n$ are freely homotopic in $ C(\overline{M}_n)$, there is an annulus $\bar E_n$ bounded by $\bar l_n$ and $\bar f_n$.\\

\indent
If $\bar y_n$ lies in $\bar d'_n$, then  $\bar l_n$ and $\bar f_n$ lie in different components of $\partial\overline{M}$. 
In particular $\bar E_n$ is an essential annulus for any $n$. 
Furthermore, we have 
%% $i(\bar\lambda_n,\partial \bar A_n)\leq i(\bar l_n,\bar\lambda_n)+i(\bar f_n, \bar\lambda_n)\longrightarrow 0$.
$i(\bar\lambda_n,\partial \bar E_n)\leq i(\bar l_n,\bar\lambda_n)+i(\bar f_n, \bar\lambda_n)\longrightarrow 0$.

Consider the projection $E_n$ of $\bar E_n$ to $C(M_n)$. Although $E_n$ may not be embedded, it follows from the Annulus Theorem \cite{foret} that any neighbourhood of $E_n$ contains an embedded annulus which we still denote by $E_n$. We have then $i(\lambda_n,\partial E_n)\longrightarrow 0$.

Thus we have proved:

\begin{claim}		\label{baranneau}
Let $M_n$ be a converging sequence of quasifuchsian manifolds with particles with the same topological type and converging angles.
Consider a geodesic arc $k_n\subset C(M_n)$ joining the two components of $\partial C(M_n)$ such that $\ell_{m_n}(k_n)\longrightarrow 0$. Then there is a sequence of essential annuli $E_n$ such that, up to extracting a subsequence, $i(\lambda_n,\partial E_n)\longrightarrow 0$.\hfill $\Box$
\end{claim}

\indent
If $\bar y_n$ lies in $\bar d_n$ then $\bar l_n$ and $\bar f_n$ are homotopic on $\partial C(\overline{M}_n)$. We are going to show that in this case $\lambda_n$ tends to have a leaf with a weight greater than or equal to $\pi$. 

\begin{claim}   \label{segmentpointu}
Let $M_n$ be a sequence of quasifuchsian manifolds with converging angles and let $\overline{M}_n$ be a negatively curved branched cover of $M_n$ such that the topological type of $\overline{M}_n$ does not depend on $n$. Consider a $\bar g_n$-geodesic arc $\bar k_n\subset C(\overline{M}_n)$ such that $\ell_{\bar g_n}(\bar k_n)\longrightarrow 0$ and that there is a $\bar m_n$-geodesic arc $\bar\kappa_n\in \partial C(\overline{M}_n)$ that is homotopic to $\bar k_n$ relative to its boundary and that satisfies $\ell_{\bar m_n}(\bar \kappa_n)\longrightarrow\infty$. Then %% we have 
$\liminf i(\bar \kappa_n, \bar\lambda_n)\geq\pi$.
\end{claim}

\begin{proof} 
The curve $\bar k_n\cup\bar \kappa_n$ is a skew polygon (up to approximating $\bar\kappa_n$ by piecewise geodesic segments) and bounds a disc in $ C(\overline{M}_n)$. Consider the geodesic cone $\bar D_n$ from $\bar x_n$ to $\bar k_n\cup\bar \kappa_n$. As in the proof of Lemma \ref{lm:annulus}, the induced metric on $\bar D_n$ is a hyperbolic metric with cone singularities with cone angles of at least $2\pi$. Since $\bar k_n$ is short, the local support planes at the endpoints $\bar x_n$ and $\bar y_n$ of $\bar k_n$ are close to each other (compare with the proof of Claim \ref{anneauplat}). It follows that the sum of the internal angles of $\bar D_n$ at $\bar x_n$ and $\bar y_n$ is close to being greater than $\pi$, namely there is $\eps_n\longrightarrow 0$ such that the sum of these $2$ angles is greater than $\pi-\eps_n$.  
Now the Gauss-Bonnet formula shows that $\liminf i(\bar \kappa_n,\bar \lambda_n)\geq\pi$.
\end{proof}

Using this Claim we will now show that under the right hypothesis $\lambda_n$ tends to have a leaf with a weight greater than or equal to $\pi$.

\begin{claim}    \label{leafpi}
Let $M_n$ be a sequence of quasifuchsian manifolds with particles with the same topological type $S\times I$ and converging angles. Let $d\subset S$ be a simple closed curve and consider its geodesic representative $d_n$ on one component of $\partial C(M_n)$. Consider an arc $\kappa_n\subset d_n$ and denote by $k_n$ the geodesic arc in $M_n$ in the homotopy class of $\kappa_n$ relative to its boundary. If $\ell_{g_n}(k_n)\longrightarrow 0$, $\ell_{m_n}(\kappa_n)\longrightarrow\infty$ and $\lambda_n$ converges, its limit $\lambda_\infty$ has a leaf with a weight greater than or equal to $\pi$.
\end{claim}

\begin{proof}
Let $\overline{M}_n$ be a negatively curved branched cover of $M_n$ whose topological does not depend on $n$. Let $\bar\kappa_n$ and $\bar k_n$ be lifts of $\kappa_n$ and $k_n$ respectively. Let $\bar x_n$ and $\bar y_n$ be the endpoints of $\bar\kappa_n$ and let $\bar l_n\subset \partial C(\overline{M}_n)$ be a shortest geodesic loop based at $\bar x_n$. Let $\bar f_n\subset\partial C(M_n)$ be the geodesic loop based at $\bar y_n$ that is homotopic to $\bar l_n$ on $\partial C(\overline{M}_n)$. Let $\tilde S_n$ be the universal cover of the connected component of $\partial C(\overline{M}_n)$ containing $\bar x_n$ endowed with the induced metric. Pick a connected component $\tilde l_n\subset \tilde S_n$ of the preimage of $\bar l_n$ under the covering projection. This broken geodesic $\tilde l_n$ is invariant under a primitive covering transformation $\gamma_n$ and we denote by $\tilde f_n$ the component of the preimage of $\bar f_n$ that is also invariant under $\gamma_n$. The line $\tilde l_n$ and $\tilde f_n$ are disjoint and bound an infinite band $\tilde B_n$, they are connected by a lift $\tilde\kappa_n$ of $\bar\kappa_n$ and by its translates by $\gamma_n^k,k\in\Z$. Pick a simple closed geodesic $\bar e_n\subset\partial C(\overline{M}_n)$ and let $\tilde e_n$ be a lift of $\bar e_n$ to $\tilde S_n$.  It is easy to check that a component  $\tilde a_n$ of $\tilde e_n\cap \tilde B_n$ which is an arc connecting $\tilde l_n$ to $\tilde f_n$ satisfies $i(\tilde a_n,\tilde \lambda_n)\geq i(\tilde \lambda_n,\tilde \kappa_n)-\sharp\{\tilde a_n\cap\bigcup_k g_n^k\tilde\kappa_n\}(i(\bar l_n,\bar\lambda_n)+i(\bar f_n,\bar\lambda_n))$ where $\tilde\lambda_n$ is the preimage of $\bar\lambda_n$ under the covering projection. By Claim \ref{anneauplat}, $i(\bar l_n,\bar\lambda_n)\longrightarrow 0$ and $i(\bar f_n,\bar\lambda_n)\longrightarrow 0$ and by Claim \ref{segmentpointu}, $\liminf i(\bar\kappa_n,\bar\lambda_n)\geq\pi$. Notice that  $\sharp\{\tilde a_n\cap\bigcup_k g_n^k\tilde\kappa_n\}$ is bounded except if $\bar e_n$ spirals more and more toward $\bar c_n$. If, for example, we assume that $\bar e_n$ and $\bar c_n$ 
converge to intersecting geodesic laminations, %% lamination, 
$\bar e_n$ does not spiral toward $\bar c_n$ and we find that %% get 
$\liminf i(\bar e_n,\bar\lambda_n)\geq i(\bar e_n,\bar c_n)\pi$.

We will now use this inequality to conclude that, when $\bar\lambda_n$ converge, its limit $\bar\lambda_\infty$ has a leaf with a weight greater than or equal to $\pi$.

Let us notice that, up to extracting a subsequence, the homotopy class of $\bar c_n$ does not depend on $n$. Otherwise there is a simple closed curve $\bar e\subset\bar S$ such that $i(\bar e,\bar c_n)\longrightarrow\infty$. To see that, extract a subsequence such that $\bar c_n$ converges in the Hausdorff topology, pick a simple closed curve $\bar e$ that intersects this limit transversally and apply the inequality above to $\{\bar e_n\}=\{\bar e\}$. But $i(\bar e,\bar c_n)\longrightarrow\infty$ would contradict the assumption that $\bar\lambda_n$ converge.

Let $\bar c\subset S$ be a simple closed curve in the homotopy class defined by $\bar c_n$. By the inequality above, we have $\liminf i(\bar e,\bar\lambda_n)\geq i(\bar e,\bar c)\pi$ for any simple closed curve $\bar e$. It follows easily that $\bar c$ is a leaf of $\bar\lambda_n$ with a weight greater than or equal to $\pi$. Taking the quotient, we conclude that $\lambda_\infty$ has a a leaf with a weight greater than or equal to $\pi$.
\end{proof}

\indent
It is now easy to conclude the proof of Lemma \ref{annulusorpi}. 
Under the assumptions of Lemma \ref{annulusorpi}, namely when there is a simple closed curve $d$ such that $\ell_{m_n}(d_n)\longrightarrow \infty$, it follows from Claims \ref{petitarc}, \ref{baranneau} and \ref{leafpi} that either there is a sequence of essential annuli $E_n\subset C(M_n)$ such that $i(\lambda_n,\partial A_n)\longrightarrow 0$ or $\lambda_\infty$ contains a leaf with a weight equal to at least $\pi$.
\end{proof}

We can now deduce from Lemma \ref{annulusorpi} that under the assumptions of 
Lemma \ref{lm:compact}, the sequence of induced metrics $(m_n)_{n\ni \N}$ on $\partial C(M_n)$ is bounded.

\begin{lemma}
%%Under the assumptions of Lemma \ref{lm:compact}, the induced metric $m_n$ on $\partial C(M_n)$ is bounded.
Let $M_n$ be a sequence of quasifuchsian manifolds with particles with the same topological type and converging angles. Let $\lambda_n$ be the measured bending 
laminations on the boundary of the convex core of $M_n$ and suppose
that $\lambda_n\rightarrow \lambda_\infty$. Let $\lambda^\pm$ be the respective restrictions %% intersections 
of $\lambda_\infty$ to %% with 
the two components of the boundary $\partial C(M_n)$ of the convex core of $M_n$. Suppose that 
\begin{itemize}
\item $\lambda_-$ and $\lambda_+$ fill $S$,
\item each closed curve in the support of $\lambda_-$ (resp. $\lambda_+$)
has weight less than $\pi$.
\end{itemize}
Then the sequence of induced metrics $(m_n)_{n\in \N}$ 
on $\partial C(M_n)$ is bounded.
\end{lemma}

\begin{proof}
If $(m_n)$ is unbounded, then there is a simple closed curve $d\subset S$ with geodesic representative $d_n\subset\partial C(M_n)$ such that $l_{m_n}(d_n)$ is unbounded. By Lemma \ref{annulusorpi} and the assumptions on $\lambda$, there is a sequence $E_n$ of essential annuli such that $i(\lambda_n,\partial E_n)\longrightarrow 0$. Such a sequence of annuli contradicts the assumption that $\lambda_-$ and $\lambda_+$ fill $S$.
\end{proof}

\subsection{Convergence of convex cores}

The last step in the proof of Lemma \ref{lm:compact} is to show that, under the assumption 
that the sequence of metrics on the boundary are bounded, 
a subsequence of convex cores converges for the bilipschitz topology. Before starting the proof, 
we will discuss the %% ajoute "the"
Margulis Lemma for quasifuchsian manifolds with particles. Let us first review the %% ajoute "the"
Margulis Lemma for manifolds with variable curvature.

\begin{thm}[Margulis Lemma]	\label{margulis variable}
Given $n\in\N$ there are constant $\mu=\mu(n)>0$ and $I(n) \in\N$ with the following property. %% : 
Let $X$ be %% ba 
an $n$-dimensional Hadamard manifold which satisfies %% statisfies 
the curvature condition $-1\leq K\leq 0$ and let $\Gamma$ be a discrete group of isometries acting on $X$. For $x\in X$ let $\Gamma_\mu (x)=\{\gamma\in\Gamma | d_\gamma(x)\leq\mu\}$ be the subgroup generated by the elements $\gamma$ with $d_\gamma(x)\leq\mu$. Then $\Gamma_\mu (x)$ is almost nilpotent, thus it contains a nilpotent subgroup of finite index. The index is bounded in $I(n)$.
\end{thm}

This statement is taken from \cite[\textsection 8.3]{ballmann-gromov-schroeder}. Since we are considering a manifold $M$ homeorphic to $S\times I$, an almost nilpotent subgroup of $\pi_1(M)$ is cyclic.

Theorem \ref{margulis variable} does not hold for hyperbolic manifolds with cone singularities since the curvature is not defined at the singularities. On the other it is not hard to replace the metric in a neighborhood of the singular locus with a Riemannian metric. Furthermore, if the cone angle is at least $2\pi$, one can choose the Riemannian metric so that it has negative curvature with a lower bound depending on the cone angles and the choice of the neighborhood of the singularities. Now we consider a quasifuchsian manifold with particle $M$ and a negatively curved branched cover $\overline{M}$ of $M$. By Lemma \ref{lm:dist3}, there are $R,\eps>0$ such that any closed curve with length at most $\eps$ is at distance at least $R$ from the singularities. We replace the $R$-neighborhood of the singularities with a smooth Riemannian metric and apply Lemma \ref{margulis variable} to the resulting manifold. Notice that the lower bound on the curvature of the Riemannian metric thus obtained will depend on $R$ and the cone angle. Thus we get $\eps$ depending on $R$ and the cone angles so that for a given point $x$ in the universal cover $\widetilde{M}$ of $\overline{M}$ the subgroup of $\pi_1(\overline{M})$ generated by the set $\{\gamma\in\pi_1(\overline{M})|d(x,\gamma x)\leq\eps\}$ is cyclic. Since $\overline{M}$ is a finite branched cover, we have a similar statement for $M$, replacing $\eps$ with $\eps/p$ where $p$ is the index of the cover which depends on the cone angles of the singularities of $M$. It follows that we have a Margulis decomposition for quasifuchsian manifolds with particles:

\begin{lemma}	\label{lm:margulis}
Let $M_n$ be a sequence of quasifuchsian manifolds with particles with the same topological type and converging angles. There is $\eps$ depending on the limit angles such that, for $n$ large enough, each component of the $\eps$-thin part of $M_n$ is a neighborhood of a closed geodesic.
\end{lemma}

Notice that the Margulis tubes %% tube 
we obtain here are disjoint from the singularities so they are isometric to regular neighborhoods %% neighborhood 
of geodesics in hyperbolic $3$-manifolds. We can now discuss the convergence of quasifuchsian manifolds with particles.

\begin{lemma} \label{lm:conv}
Let $M_n$ be a sequence of quasifuchsian manifolds with particles with the same topological type and converging angles. 
Suppose that the sequence $(m_n)_{n\in \N}$ %% $\{m_n\}$ 
(the induced metrics on the boundary of the convex cores) converges. 
Then, after taking a subsequence, $(M_n)_{n\in \N}$ converges to a quasifuchsian manifold with particles with the same topological type as $M_n$.
\end{lemma}

%% dans cette partie j'ai remplace m_n par g_n pour garder la coherence des notations.

\begin{proof}
First notice that since the cone angles are less than $\pi$, by Lemma \ref{lm:dist3} there is a positive lower bound 
for the distance between two components of the singularity locus. 
Consider a point $x_n\in C(M_n)$, extract a subsequence such that the sequence $(x_n, M_n)$ 
converges in the Gromov-Hausdorff topology (such a subsequence always exists). 
By \cite[Proposition 3.2.6]{boileau-porti}, the limit $(x_\infty, M_\infty)$ 
is a hyperbolic manifold with cone singularities. 
By \cite[Proposition 3.3.1]{boileau-porti}, the sequence $(x_n, M_n)$ converges to  
$(x_\infty, M_\infty)$ in the bilipschitz topology. 
%It follows easily that $ M_\infty$ is a convex manifold with cone singularities. The bilipschitz convergences implies that the cone singularities of the limit manifold are orthogonal to the boundary.

It remains to show that $M_\infty$ has the same topological type as $M_n$ and that its metric is convex co-compact. To do that we will show that the diameter of $ C(M_n)$ is uniformly bounded. It will follow that $C(M_n)$ converges to a convex set with the same topological type.

\begin{lemma}   \label{diamborne}
Let $(M_n)_{n\in \N}$ be a sequence of quasifuchsian manifolds with particles with the same topological type and converging angles. 
Suppose that the sequence of induced metrics $(m_n)_{n\in\N}$ on the boundary of the convex cores converges. 
Then the diameter of $C(M_n)$ is uniformely bounded.
\end{lemma}

\begin{proof}
Consider a negatively curved ramified cover $\overline{M}_n$ of $M_n$ whose topological type does not eventually depend on $n$. 
It follows from the %% ajoute
Margulis Lemma that a very short geodesic in $\overline{M}_n$ 
lies in a very deep embedded tube. 
Using this %% that 
observation we will show that there is a uniform lower bound on the length
of any fixed curve in $ C(\overline{M}_n)$.

\begin{claim} \label{paspetit}
Let $M_n$ be a sequence of quasifuchsian manifolds with particles and let $\overline{M}_n$ be 
a negatively curved ramified cover of $M_n$ whose topological type does not depend on $n$. 
Suppose that the sequence $(\bar m_n)_{n\in \N}$ %% $\{\bar m_n\}$ 
converges. Let $\bar c\subset \bar S$ be a simple closed curve. 
Then there is $Q>0$ such that if $\bar c_n\subset C(\overline{M}_n)$ 
denotes the geodesic representative of $\bar c$, $l_{\bar g_n}(\bar c_n)\geq Q$ for any $n\in \N$. 
\end{claim}

\begin{proof}
Assume the contrary, that is (after extracting a subsequence), $\lim l_{\bar g_n}(\bar c_n)=0$. 
Then $\bar c_n$ is the core of a deep Margulis tube $\overline{T}_n$. 
Notice that since $(\bar m_n^+,\bar m_n^-)_{n\in \N}$ converges, there is no short curve in $\partial C(\overline{M}_n)$. More precisely, there is a uniform lower bound on the length of simple closed geodesics on $\partial C(\overline{M}_n)$.  Since the induced metric on $\partial C(\overline{M}_n)$ is negatively curved, it can have a large diameter only if it contains a short curve. Thus the uniform lower bound on the length of simple closed geodesics on $\partial C(\overline{M}_n)$ provides us with a bound on the diameter of each component of  $\partial C(\overline{M}_n)$. It follows that $\partial C(\overline{M}_n)$ does not go too deep into a Margulis tube (compare with  \cite[Lemma 6.3]{minsky}): Let $\eps_0$ be a Margulis constant for the sequence $M_n$ as provided by Lemma \ref{lm:margulis}, namely the $\eps_0$-thin part $M_n^{<\eps_0}$ of $M_n$ is a union of Margulis tubes for $n$ large enough. By \cite{meyerhoff} and \cite{brooks-matelski} (see also \cite[Lemma 6.1]{minsky}) given $\eps$ small enough, the diameter of  $M_n^{<\eps_0}-M_n^{<\eps}$ is large. In particular, if a component of $\partial C(\overline{M}_n)$ intersects $M_n^{<\eps}$ for a small $\eps$, it has a large diameter. Hence the bound on the diameter of each component of  $\partial C(\overline{M}_n)$ provides us with a constant $\eps$ (depending on the sequence $M_n$) so that $\partial C(\overline{M}_n)$ is disjoint from the thin part $M_n^{<\eps}$.

If we take  $\overline{T}_n$ to be an $\eps$-Margulis tube, 
we get that $\overline{T}_n$ lies entirely in the interior of $ C(\overline{M}_n)$. 
Consider a simple closed curve $\bar d\subset \bar S$ that intersects $\bar c$ essentially. By Lemma \ref{lm:annulus}, there is an essential annulus $\overline{A}_n\subset C(\overline{M}_n)$ 
which is in the homotopy class defined by $\bar d\times I$ such that the area of $\overline{A}_n$ is at most $\ell_{\bar m_n^+}(\bar d_+)+\ell_{\bar m_n^-}(\bar d_-)$. In particular, since the sequence $(\bar m_n)_{n\in \N}$ %% $\{\bar m_n\}$ 
converges, the area of $\overline{A}_n$ is bounded. On the other hand, since $\bar d$ intersects $\bar c$ essentially, $\overline{A}_n$ intersects $\bar c_n$ essentially. In particular, $\overline{A}_n$ intersects $\overline{T}_n$ along a disc $\overline{D}_n$. When the length of $\bar c_n$ tends to $0$, $d(\bar c_n,\partial \overline{T}_n)\longrightarrow\infty$ (see \cite{meyerhoff}, \cite{brooks-matelski} and \cite[Lemma 6.1]{minsky}). It follows that the diameter of $\overline{D}_n$, and hence its area, tends to $\infty$ when the length of $\bar c_n$ tends to $0$. Thus an upper bound for the area of $\overline{A}_n\supset\overline{D}_n$ yields a lower bound for the length of $\bar c_n$. This concludes the proof of Claim \ref{paspetit}.
\end{proof}

Consider now two simple closed curves $\bar c,\bar d\subset\bar S$ such that the components of $\bar S\setminus (\bar c\cup\bar d)$ are discs. 
Two such curves are said to fill the surface $\bar S$. 
Consider essential annuli $\overline{A}_n$ and $\overline{B}_n$ in $ C(\overline{M}_n)$ in the homotopy classes defined by $\bar c$ and $\bar d$, constructed as in Lemma \ref{lm:annulus}. 
In particular $\overline{A}_n$ and $\overline{B}_n$ have bounded area. Since $\overline{A}_n$ and $\overline{B}_n$ have bounded areas and negative curvature, the only way for them to have a large diameter is to have a very short core curve. 
But this would contradict Claim \ref{paspetit}. Thus we can conclude that $\overline{A}_n$ and $\overline{B}_n$ have uniformly bounded diameters.

Let ${\cal B}_{1\leq k\leq p}$ be the closure of the components of 
$\bar S\times I\setminus (\bar c\times I\cup\bar d\times I)$. 
Our manifold $\bar N=\bar S\times I$ is the union of the ${\cal B}_k$ and the ${\cal B}_k$ are all balls. Define a surjective map $f_n:\bar S\times I\rightarrow  C(\overline{M}_n)$ that maps $\bar c\times I$ and $\bar d\times I$ to $\overline{A}_n$ and $\overline{B}_n$ respectively and such that the restriction of $f_n$ to each ${\cal B}_k$ is an immersion. For each $k$, the image of $\partial{\cal B}_k$ lies in $\overline{A}_n\cup \overline{B}_n\cup\partial C(M_n)$. Since $\overline{A}_n$ and $\overline{B}_n$ have bounded diameters and since the induced metric on $\partial C(\overline{M}_n)$ is bounded, the diameter of $f_n(\partial{\cal B}_k)$ is bounded for any $k$. It follows that $f_n({\cal B}_k)$ has a bounded diameter for any $k$. Since $f_n$ is surjective, this implies that $ C(\overline{M}_n)$ has a uniformly bounded diameter. Since the index of the cover $\overline{M}_n\rightarrow M_n$ does not depend on $n$, $C(M_n)$ has %% have 
a bounded diameter.
\end{proof}

\indent
It remains to show that the convex core of $M_\infty$ is compact and homeomorphic to $S\times I$ (when $M_n$ is homeomorphic to $S\times\R$). Once again we will use the %% ajoute
negatively curved ramified cover $\overline{M}_n$. Since $(\bar x_n,\overline{M}_n)$ converge to $(\bar x_\infty,\overline{M}_\infty)$, there is $R_n\longrightarrow\infty$ and a sequence of bilipschitz map $\phi_n:B(\bar x_n,R_n)\rightarrow B(\bar x_\infty,R_n)$ such that the bilipschitz constants tend to $1$. By Lemma \ref{diamborne}, for $n$ large enough, $C(\overline{M}_n)\subset B(x_n,R_n)$. Given a geodesic segment $\gamma_n\subset C(\overline{M}_n)$, $\phi_n(\gamma_n)$ almost realizes the distance between its endpoints. Since $\overline{M}_n$ is a hyperbolic manifold with cone singularities, $\phi_n(\gamma_n)$ is very close to the geodesic segment joining its endpoints. It follows that, for $n$ large enough, the convex hull of $\phi_n(C(\overline{M}_n))$ lies in a small neighborhood ${\cal V}_n(\phi_n(C(\overline{M}_n)))$ of $\phi_n(C(\overline{M}_n))$. This convex hull has to contain $C(\overline{M}_\infty)$ since it is the smallest convex set. Thus we have $C(\overline{M}_\infty)\subset {\cal V}_n(\phi_n(C(\overline{M}_n)))$. It follows that $C(\overline{M}_\infty)$ is compact. Furthermore, since the induced metric on $\partial C(\overline{M}_n)$ is bounded, ${\cal V}_n(\phi_n(C(\overline{M}_n)))$ is homeomorphic to $\bar S\times I$ for $n$ large enough. It follows that $C(\overline{M}_\infty)$ is homeomorphic to $\bar S\times I$. Thus we have proved that $M_\infty$ is a quasifuchsian manifold with cone singularities with the same topological type as $M_n$.
\end{proof}

In contrast to the other results of this section, we do not need to assume that the quasifuchsian manifolds under consideration in Lemma \ref{lm:conv} are not fuchsian.

\subsection{The bending lamination of the convex core}

To finish the proof of Lemma \ref{lm:compact} we only have to check that
the induced bending lamination on the boundary of the convex core of the
limit manifold is the limit of the bending laminations. We can state
the result as follows.

\begin{lemma} \label{lm:limit-bending}
Let $N=S\times \R$, let $x_1,\cdots, x_{n_0}$ be distincts points on $S$, and let  
$\kappa_i=\{ x_i\}\times \R, 1\leq i \leq k$. 
Let $(g_n)_{n\in \N}$ be a sequence of quasifuchsian metrics on $N$ with particles
of angles $\theta^i_n$ along $\kappa_i$, $1\leq i\leq k$. 
Let $\lambda_n$, resp. $m_n$, be the measured bending 
laminations, resp. the induced metric, on the boundary of the convex core of $(N,g_n)$.
Suppose that $(g_n)$ converges in bilipschitz topology towards a quasifuchsian
metric with particles $g$ on $N$, with cone angles $\theta^i\in (0,\pi)$ along $\kappa_i$.
Then $(m_n)_{n\in \N}$ converges to the induced metric $m$ on the boundary of the convex core
of $(N,g)$, while $(\lambda_n)_{n\in \N}$ converges to the measured bending lamination
$\lambda$ of the boundary of the convex core of $(N,g)$. 
\end{lemma}

\begin{proof}
Set $M=(N,g)$ and $M_n=(N,g_n)$ and denote by $C(M)$ the convex core of $M$ and by $C(M_n)$ the convex core of $M_n$.
We consider as above the finite cover $\bar N$ of $N$ ramified along the cone
singularities, chosen so that all cone angles in $\bar N$ have angle larger than
$2\pi$. This is useful below since we will use negative curvature arguments, in
particular the existence of a geodesic segment in a homotopy class with fixed
endpoints. Clearly it is sufficient to prove the lemma for $\bar N$, where the
``convex core'' considered is $ C(\overline{M}_n)$, the lift to $\bar N$ of $ C(M_n)$,
since once the result is obtained in $\bar N$, we can take the quotient by
the group of deck transformations to obtain the result on $N$.

Let $(\gamma_n)_{n\in \N}$ be a sequence of segments in $\bar N$, with $\gamma_n$
geodesic for $\bar g_n$ for all $n\in \N$. 
Suppose that $(\gamma_n)_{n\in \N}$ converges to a segment $\gamma$.
We know that $\bar g_n\rightarrow \bar g$ in the bilipschitz topology and, in hyperbolic geometry,
any segment which is close to realizing the distance between its endpoints is close 
to a geodesic segment. So $\gamma$ is geodesic for $\bar g$. Conversely, any geodesic
segment for $\bar g$ is a Hausdorff %% Haudorff 
limit of geodesic segments for the $g_n$. The same holds for
closed geodesics.

Let $(\Omega_n)_{n\in \N}$ be a sequence of compact subsets of $\bar N$ such that, for all 
$n\in \N$, $\Omega_n$ is convex for $\bar g_n$. Suppose that $\Omega_n\rightarrow \Omega$
in the Hausdorff topology. 
The definition of a convex subset and the previous paragraph show that $\Omega$ 
is convex, since any geodesic segment $\gamma$ in $(\bar N,\bar g)$ with endpoints in 
the interior of $\Omega$ is
the limit of a sequence of geodesic segments $\gamma_n$, with $\gamma_n$ geodesic for
$\bar g_n$. Since $\gamma_n$ has endpoints in $\Omega_n$ (for $n$ large enough) and 
$\Omega_n$ is convex for $\bar g_n$, $\gamma_n\subset\Omega_n$, and therefore 
$\bar \gamma\subset \Omega$, and $\Omega$ is convex for $\bar g$. 
Conversely, a similar argument shows that any compact convex subset for $\bar g$ is the 
Hausdorff limit of a sequence of compact convex subsets of the metrics $\bar g_n$.

For all $n$, $ C(\overline{M}_n)$ contains all closed geodesics in $(\bar N, \bar g_n)$. Given a 
non-trivial homotopy class $\alpha$ in $\bar N_r$ (the complement of the singular curves 
in $\bar N$), it is realized for each $n\in \N$ by a (unique) closed geodesic $\gamma_n$ in 
$(\bar N,\bar g_n)$, and the sequence $(\gamma_n)_{n\in \N}$ converges to the closed geodesic $\gamma$
which realizes $\alpha$ in $(\bar N,\bar g)$. For each $n\in \N$, $\gamma_n\subset  C(\overline{M}_n)$.
Moreover we have seen that the diameter of the $ C(\overline{M}_n)$ is bounded. It follows that
$( C(\overline{M}_n))_{n\in \N}$ converges -- after extracting a subsequence -- to a limit 
subset $C'$ which contains all closed geodesics in $(\bar N,\bar g)$. 

Since $C'$ is the limit of a sequence of convex subset of the $(\bar N,\bar g_n)$, it is
convex. Moreover
if $\Omega\subset C'$ is convex, then it is the limit of a sequence of convex 
subsets $\Omega_n\subset  C(\overline{M}_n)$. But then $(\Omega_n\cap  C(\overline{M}_n))_{n\in \N}$ 
is a sequence of convex subsets converging to $\Omega$. Because the $ C(\overline{M}_n)$
are minimal convex subsets, $\Omega_n\cap  C(\overline{M}_n)= C(\overline{M}_n)$ for all $n$, so 
that $\Omega=C'$. So $C'= C(\overline{M})$. This shows that $ C(\overline{M})$ is the Hausdorff
limit of the $ C(\overline{M}_n)$. 

Note that it is not clear at this point whether $\partial C(\overline{M}_n)\rightarrow \partial  C(\overline{M})$
in the $C^1$ topology. 
However, a general fact is that, if $\phi:S\rightarrow H^3$ is a smooth embedding of a
surface, and if $(\phi_n)_{n\in \N}$ is a sequence of Lipschitz embeddings of
$S$ in $H^3$ which converges to $\phi$ in the $C^0$ topology, then the distance
$d_n$ induced on $S$  by the $\phi_n$ are larger, in the limit, than
the distance $d$ induced by $\phi$:
\begin{equation}\label{eq:limsup}
  \forall x,y\in S, \limsup_{n\rightarrow \infty} d_n(x,y)\geq
d(x,y)~. 
\end{equation}
The same holds if $\phi$ is Lipschitz with locally convex image rather than smooth, 
see \cite{AZ}. Moreover, in case of equality in Equation (\ref{eq:limsup}) and if the 
$\phi_n$ also have locally convex images, then
the convergence of $(\phi_n)$ to $\phi$ is stronger, in the sense that the tangent
plane to $\phi_n(S)$ almost everywhere converges to the tangent plane to $\phi(S)$.

Coming back to $\partial C(\overline{M}_n)$, 
the $C^0$ convergence towards $\partial  C(\overline{M})$
(together with the bilipschitz convergence
of $\bar g_n$ to $\bar g$) is sufficient to insure that the metric $\bar m_n$ on $\dr  C(M_n)$
is {\it larger} in the limit than the metric $\bar m$ induced by $\bar g$ on $\dr  C(\overline{M})$. In other terms,
if $x,y\in \dr  C(\overline{M})$ and $x_n,y_n\in \dr  C(\overline{M}_n)$ are such that $\lim x_n=x, \lim y_n=y$,
then there exists for each $\epsilon>0$ some $N_0\in \N$ such that, for 
all $n\geq N_0$, 
$$ d_{\bar m_n}(x_n,y_n)\geq (1-\epsilon)d_{\bar m}(x,y)~. $$
It follows that the lengths of the closed geodesics in $S$ for $\bar m_n$ are 
bounded from below by $(1-\epsilon)$ times their lengths for $\bar m$. 

Since the metrics $\bar m_n$ and $\bar m$ are hyperbolic metrics with cone singularities of 
fixed angles, this shows, using standard arguments based for instance on pants
decompositions, that $\bar m_n\rightarrow \bar m$ (see Section \ref{ap:A2}, or
\cite{dryden-parlier}).

It then follows that as $n\rightarrow \infty$, $\bar m_n$ is bounded from above by
$(1+\epsilon)\bar m$: under the same hypothesis as above, 
\begin{equation}\label{eq:sup}
 d_{\bar m_n}(x_n,y_n)\leq (1+\epsilon)d_{\bar m}(x,y)~.
\end{equation}

Since $\bar m_n\rightarrow \bar m$ in this sense of Equation (\ref{eq:sup}), 
and moreover $\dr C(\overline{M}_n)\rightarrow \dr  C(\overline{M})$
in the $C^0$ topology, it follows that the convergence is actually stronger,
and the tangent plane to $\dr C(\overline{M}_n)$ converges almost everywhere to the tangent plane
to $\dr C(\overline{M})$ (both exist almost everywhere by convexity). 
This implies that the measured laminations $\lambda_n$ 
of $\dr C_n$ converge to $\lambda$.
\end{proof}

%% file: conebend4.tex
\section{Prescribing the measured bending lamination on the boundary of the convex core}
\label{se:rational}

The goal of this section is to prove Theorem \ref{tm:rational} and then
Theorem \ref{tm:general}. The proof of Theorem \ref{tm:rational} is
largely based on a well-known doubling argument already used for non-singular 
manifolds, which reduces the infinitesimal rigidity with respect to the
measured lamination (when the support of the lamination is along closed
curves) to a rigidity statement proved by Hodgson and Kerckhoff \cite{HK} for 
hyperbolic cone-manifolds. 

Theorem \ref{tm:general} is then a consequence, using the compactness
statement proved in section \ref{se:compact-lamination}.

\subsection{A doubling argument}

Let $M$ be convex co-compact manifold with particles, and let $C(M)$ be
its convex core. Suppose that the support of the measured bending lamination of 
$C(M)$ is a disjoint union of closed curves.

\begin{df}
The {\bf doubled convex core} of $M$ is the 3-dimensional hyperbolic
manifold with cone singularities $DC(M)$ obtained by gluing two copies of $C(M)$
isometrically using the identification of their boundaries.
\end{df}

We have seen that the singular locus of $M$ does not intersect the support of 
the bending lamination on the boundary of the convex core --- actually it even
remains at a distance which is bounded from below by a constant depending only
on the cone angles. So the ``particles'' intersect the boundary of the convex core
inside faces, and moreover it does so orthogonally. 
It follows that the singular locus of $DC(M)$ is a disjoint
union of closed curves, which are of two types:
\begin{itemize}
\item each ``particle'' $p$ of $M$ corresponds to a cone singularity
along a closed curve in $DC(M)$, of length equal to twice the length of
the intersection of $p$ with $C(M)$,
\item each closed curve in the support of the measured bending lamination of
the boundary of $C(M)$ corresponds to a closed curve (of the same length)
in $DC(M)$. 
\end{itemize}
Still by definition, $DC(M)$ admits an isometric involution --- exchanging
the two copies of $C(M)$ which are glued to obtain $DC(M)$ --- and the 
set of fixed points of this involution is a (non connected)) closed surface $S$, which
corresponds to the boundaries of both copies of $C(M)$. This surface is orthogonal to
the singularities of the first kind, and contains the singularities of the second
kind. 

\subsection{Local deformations}

The doubling trick explained above leads directly to a rigidity statement. We consider
again a convex co-compact manifold $M$ with particles, for which the
measured bending lamination of the convex core is along closed curves $\gamma_1,\cdots, \gamma_N$,
for which the bending angles are equal to $\alpha_1,\cdots, \alpha_N\in (0,\pi)$. As in the
introduction, we call $\theta_1,\cdots,\theta_{n_0}$ the cone angles at the ``particles'', and
let $\theta=(\theta_1,\cdots,\theta_{n_0})$.

%% corrige le lemme suivant, plein de bugs.
\begin{lemma} \label{lm:local-rational}
There exists a neighborhood $U$ of $(\alpha_1,\cdots,\alpha_N)$ in $(0,\pi)^N$
and a neighborhood $V$ of the hyperbolic metric $g$ on $M$ in $\QF_{S,n_0,\theta}$
such that, if $(\alpha'_1,\cdots,\alpha'_N)\in U$, there is a unique $g'\in V$
for which the support of the measured bending lamination on $C(M)$ is
$\gamma_1\cup\cdots\cup \gamma_N$ and the bending angle on $\gamma_i$ 
is $\alpha'_i$, $1\leq i\leq N$. 
\end{lemma}

\begin{proof}
Hodgson and Kerckhoff \cite{HK}  proved a local deformation result for hyperbolic
cone-manifolds. It follows from their result that there exists a unique cone-manifold
close to $DC(M)$ with the same topology as $DC(M)$ (including the singular locus),
the same angles at the cone singularities corresponding to the particles in $M$,
and angles $2\alpha'_1,\cdots, 2\alpha'_N$ instead of $2\alpha_1,\cdots, 2\alpha_N$
at the cone singularities corresponding to the pleating lines of $C(M)$.

The uniqueness of $D'$ shows that it has the same symmetry as $DC(M)$, that is, it 
admits an isometric involution fixing a surface $S'$ isotopic to the surface $S$
fixed by the isometric involution on $DC(M)$. By an easy symmetry argument, the
cone singularities in $D'$ corresponding to the particles in $M$ still have to be
orthogonal to $S'$, while those corresponding to the pleating lines of $\dr
C(M)$ have to be contained in $S'$ (see 
\cite[Section 8]{bonahon-otal} for details on the uniqueness part of this argument;
the same argument can basically be used when particles are present). 

Therefore, $D'$ is the double of a hyperbolic manifold with convex boundary
(obtained as the metric completion of one half of the complement of $S'$
in $D'$) with cone singularities orthogonal to the boundary. 
The boundary of this manifold is convex with no extremal point,
so that it is the convex core of a quasifuchsian manifold with particles
$M'$, with the same cone angle as $M$ at the particles and such that $\dr C(M')$
is pleated along the same lines as $\dr C(M)$, but with pleating angles 
$\alpha'_1,\cdots, \alpha'_N$ instead of $\alpha_1,\cdots, \alpha_N$.

The uniqueness of such a manifold, in the neighborhood of $M$, follows from
the uniqueness of $D'$ in the neighborhood of $DC(M)$.
\end{proof}

\subsection{Proof of Theorem \ref{tm:rational}}

Let $\gamma_1,\cdots,\gamma_N$ be the curves in the support of $\lambda$,
considered as curves in $\dr N$. Following the doubling construction above,
we define a closed manifold $D(N)$ by gluing two copies of $N$ along their
boundary. $D(N)$ contains two families of curves, which we still call 
$c_1,\cdots,c_{n_0}$ (corresponding to the particles in $N$) and 
$\gamma_1,\cdots,\gamma_N$ (corresponding to the pleating lines on the boundary
of the convex core). 

Let $\theta'_1,\cdots,\theta'_{n_0}\in (0,\pi)$ and $\alpha'_1,\cdots, \alpha'_N\in (0,\pi)$
be chosen such that: 
\begin{itemize}
\item for all $i\in \{ 1,\cdots, n_0\}$, $0\leq \theta'_i\leq \theta_i$, and $\theta'_i=\pi/k_i$
for some $k_i\in \N$,
\item for all $j\in \{ 1,\cdots, N\}$, $0\leq \alpha'_j\leq \alpha_j/2$, and $\alpha'_j=\pi/2l_j$
for some $l_j\in \N$.
\end{itemize}
The Orbifold Hyperbolization Theorem for cyclic orbifolds (initially stated by Thurston, 
and proved in \cite{boileau-porti,CHK}) can be applied to show that there is a unique
hyperbolic orbifold structure on $D(N)$ with singularities of angles $\theta'_i$
on the $c_i$ and $2\alpha'_j$ on the $\gamma_j$. 

Since the $\theta_i$ are in $(0,\pi)$, the result of Kojima \cite{kojima-cone} shows
that this orbifold structure can be deformed to a unique cone-manifold structure, with
cone angles $\theta_i$ on the curves $c_i$ and $2\alpha'_j$ on the curves $\gamma_j$.

Let $(\alpha_t)_{t\in [0,1]}=(\alpha_{1,t},\cdots,\alpha_{N,t})_{t\in [0,1]}$ be the 1-parameter 
family defined by
$$ \alpha_{j,t} = (1-t)\alpha'_j + t\alpha_j, ~ 1\leq j\leq N~. $$
Then for all $j\in \{ 1,\cdots, N\}$, $\alpha_{j,0}=\alpha'_j, \alpha_{j,1}=\alpha_j$.
Let $I\in [0,1]$ be the maximal interval containing $0$ such that, for all $t\in I$:
\begin{itemize}
\item there exists a hyperbolic structure on $D(N)$ with cone singularities of angle 
$\theta_i$ on $c_i$, $1\leq i\leq n_0$, and a cone singularity of angle $2\alpha_{j,t}$
on $\gamma_j$, $1\leq j\leq N$,
\item this hyperbolic structure has an isometric involution exchanging the 
two copies of $N$ glued to obtain $D(N)$.
\end{itemize}
By construction, $I\neq\emptyset$. Lemma \ref{lm:local-rational} shows that 
$I$ is open, while Lemma \ref{lm:compact} shows that $I$ is closed. So 
$I=[0,1]$, this proves the existence part of the statement because $D(N)$ with
the hyperbolic cone-structure for $t=1$ is obtained by doubling the convex
core of a convex co-compact hyperbolic manifold with particles of angles 
$\theta_i$ and pleating angles $\alpha_i$ on the boundary, as needed.

For the uniqueness, the same deformation argument can be used to start from a 
cone-manifold structure on $D(N)$ and decrease
the angles along the curves $\gamma_j$, $1\leq j\leq N$, from $2\alpha_j$ to
$2\alpha'_j$. Lemma \ref{lm:local-rational} shows that the corresponding deformation of 
the hyperbolic cone-manifold structure exists and is unique. Since the endpoint
of the deformation is unique (by the Orbifold Hyperbolization Theorem) there can
be only one cone-manifold structure on $D(N)$ with angles $\theta_i$ on the curves
$c_i$, $1\leq i\leq n_0$, angle $2\alpha_j$ on the curve $\gamma_j$, $1\leq j\leq N$,
and the necessary symmetry property.

\subsection{Proof of Theorem \ref{tm:general}.}

Given $\lambda_-,\lambda_+\in \cML_{S,x}$ satisfying the hypothesis of
Theorem \ref{tm:general}, both are limits of a sequence of 
measured laminations $(\lambda_{-,n})_{n\in \N},(\lambda_{+,n})_{n\in \N}$ with
support along a union of closed curves, which satisfy the hypothesis of 
Theorem \ref{tm:rational}. 

For all $n$, Theorem \ref{tm:rational} shows that $\lambda_{-,n}$ and $\lambda_{+,n}$
are the upper and lower measured bending laminations of the boundary of the convex
core for a unique quasifuchsian hyperbolic structure with particles $g_n$ on 
$S\times \R$. Lemma \ref{lm:compact}, applied to this sequence of hyperbolic
structures, shows that it has a subsequence which converges to a quasifuchsian
hyperbolic structure with particles, for which the lower and upper measured 
bending laminations of the boundary of the convex core are $\lambda_-$
and $\lambda_+$, respectively.

\subsection{The conditions are necessary}

Finally we check here that the hypothesis in Theorem \ref{tm:general} are
necessary. It obviously follows that the hypothesis in Theorem \ref{tm:rational}
are also necessary.

\begin{lemma} \label{lm:n-general}
Let $M$ be a non Fuchsian quasifuchsian manifold with particles, let $\lambda$
be the measured bending lamination on the boundary of its convex
core. Then $\lambda$ satisfies the hypothesis of Theorem \ref{tm:general}.
\end{lemma}

\begin{proof}
The hypothesis that the weight of each closed curve in the support of 
$\lambda_-$ and $\lambda_+$ is less than $\pi$ is clearly a consequence of
the fact that $C(M)$ is convex and compact. 

Suppose by contradiction that $\lambda_-$ and $\lambda_+$ do not fill $S$. There exists then
a sequence $(c_n)_{n\in \N}$ of simple closed curve in $S$ such that 
$$ i(\lambda_-,c_n)+i(\lambda_+, c_n)\rightarrow 0~. $$
Let $c_n^-$ and $c_n^+$ be the geodesic representatives
of $c_n$ in the lower and upper boundary components of $C(M)$, respectively. 

Let $\bar c_n^-$ and $\bar c_n^+$ be lifts of $c_n^-$ and $c_n^+$, respectively, to $\bar M$,
corresponding to the same lift of $c$.
Lemma \ref{lm:annulus} shows that 
there exists an annulus $A_n\subset \bar M$ bounded by $\bar c_n^-$ and $\bar c_n^+$ on
which the induced metric is hyperbolic with cone points of negative singular curvature
(cone angle larger than $2\pi$). Moreover, the boundary of $A_n$ is convex (for the 
induced metric) and its total curvature goes to $0$ as $n\rightarrow \infty$.
The Gauss-Bonnet formula then implies that the area of $A_n$ 
goes to $0$ as $n\rightarrow \infty$. Since the lengths of the 
$\bar c_n^-$ and $\bar c_n^+$ are
bounded from below, this means that the distance between $\bar c_n^-$ and $\bar c_n^+$ in 
$A_n$ goes to $0$ as $n\rightarrow \infty$. Therefore, the distance 
between the upper and lower boundary of $C(M)$ is zero, a contradiction because we
have supposed that $M$ is not Fuchsian.
\end{proof}

%% file: conebend5.tex
\section{Earthquakes estimates}
\label{se:earthquake}

In this section we consider a convex co-compact manifold with 
particles $M$. The arguments in this more general case are the
same as in the specific situation of quasifuchsian manifolds
with particles. Its boundary $\dr M$ has a number of marked
points $x_1, \cdots, x_{2n_0}$ which are the endpoints of the $n_0$
``particles'', and to each is attached an angle $\theta_k\in (0,\pi)$,
$1\leq k\leq 2n_0$, which is the angle at the corresponding particle.
%% change la numerotation pour simplifier et faire apparaitre le nombre de particules.

We identify $\dr M$ with the boundary of 
its convex core (see Lemma \ref{lm:ends}). We will use the following notations. 
\begin{itemize}
\item $\lambda$ is the measured bending lamination of the boundary of the
convex core.
\item $m$ is its induced metric.
\item  $t$ is the (unique) hyperbolic metric in the conformal class at 
infinity $\tau$, with cone singularities
of prescribed angle $\theta_k$ at the marked point $x_k$.
\item $G_\lambda(m)$ is the metric obtained by grafting
the hyperbolic metric $m$ along the measured lamination $\lambda$, so
that $G_\lambda(m)$ has curvature in $[-1,0]$. If for instance $\lambda$
is rational, then $G_\lambda(m)$ is obtained by inserting a flat annulus
in $(\dr M,m)$ for each closed curve in the support of $\lambda$, see 
e.g. \cite{dumas-survey}.
\end{itemize}
 
This section contains a basic estimates relating $t$ to $m$.
It will be useful in proving
the compactness of a certain map and Theorem \ref{tm:bers}.
Its statement is based on the following extension to hyperbolic surfaces
with cone singularities of Thurston's Earthquake Theorem (as found
in \cite{thurston-earthquake,kerckhoff}).

\begin{thm}[\cite{cone}]
For any $h,h'\in \cT_{\Sigma,x,\theta}$, there is a unique measured
lamination $\nu\in \cML_{\Sigma, x}$ such that the right earthquake 
along $\nu$ sends $h$ to $h'$.
\end{thm}

The main estimate proved in this section, and the main tool for
the proof of  Theorems \ref{tm:bers} and \ref{tm:metriques}, is the following.

\begin{prop} \label{pr:length}
There exists a constant $C>0$ (depending only on the topology
of $M$) such that, if $\nu\in \cML_{\dr M,x}$ is the measured
lamination such that $t=E_r(\nu)(m)$, then the length
$L_m(\nu)$ is at most $C$. 
\end{prop}

It is proved in Section \ref{sc:length}, after some preliminary 
considerations. It is used below in Section \ref{se:induced}.

\subsection{The average curvature of geodesics}

In this part we prove a technical statement which is useful at several
points below. It is an extension to convex co-compact manifolds with
particles of a result proved earlier by Bridgeman \cite{bridgeman} for
the convex core of non-singular convex co-compact manifolds, or more
generally of pleated surfaces in $H^3$. However the argument used here
is inspired by Bonahon and Otal \cite{bonahon-otal}.

We consider a quasifuchsian manifold with particles, $M$, and
call $\theta_1,\cdots,\theta_{n_0}$ the cone angles at the particles. By
definition, $\theta_1,\cdots, \theta_{n_0}\in (0,\pi)$. Here $S$ is one
of the connected components of $\partial C(M)$.

\begin{prop} \label{pr:average}
There exists a constant $C_0>0$ such that, if $\gamma$ is a geodesic
segment on $S$ transverse to $\lambda$, $i(\gamma,\lambda)\leq 
C_0(l_m(\gamma)+1)$.
\end{prop}

Note that $C_0$ depends on the $\theta_i$ (at least the argument we use
here does depend on the maximum of the $\theta_i$) but not otherwise
on $M$. 

Proposition \ref{pr:average} will follow from the following lemma. We use
here the constant $\epsilon_0$ coming from Lemma \ref{lm:dist2} and Lemma
\ref{lm:dist3}.

\begin{lemma}\label{lm:max}
There exists $\lambda_1>0$ such that if $\gamma$ is at distance at
least $\epsilon_0/2$ from the intersection of $S$ with the singular set of 
$M$ and if $l_m(\gamma)\leq \epsilon_0/4$, then $i(\lambda,\gamma)\leq \lambda_1$. 
\end{lemma}

The proof of this lemma is based on some intermediate steps. The first
is a consequence of Lemma \ref{lm:dist2} and Lemma \ref{lm:dist3}.

\begin{claim}
There exists $\rho_0>0$ such that any point in $S$ at distance
at least $\epsilon_0$ in $S$ from the singular points of $S$ is also at distance
at least $\rho_0$ in $M$ from the singular set of $M$. 
\end{claim}

\begin{proof}
Let $x\in S$ which is at distance at least $\epsilon_0$ from the singular
points of $S$, suppose that it is at distance strictly less than $\epsilon_0$
from a particle $p$. Let $y$ be a point in $p$ which is closest from $x$, and 
let $D$ be the totally geodesic disk of radius $\epsilon_0$ orthogonal to 
$p$ at $y$. This disk does not encounter any other particle by 
Lemma \ref{lm:dist3}. Moreover $x\in D$ because $y$ is at minimal distance 
from $x$ among the points of $p$. We can therefore apply the second point in
Lemma \ref{lm:dist2} to $D$, with $\Omega$ equal to the intersection of 
$D$ with the convex core of $M$. The result follows.
\end{proof}

\begin{cor}
Let $x\in S$ be contained in the support of the bending lamination $\lambda$, 
and let $D'$ be the totally geodesic disk of radius $\rho_0$ in $M$ orthogonal to $\lambda$ 
at $x$. Then $D'$ does not intersect the singular set of $M$.
\end{cor}

After taking $\rho_0$ smaller if necessary, we have another simple statement
which will be necessary below. 

\begin{claim} \label{cl:angle}
Let $y\in \lambda$ be a point in the connected 
component of $x$ in the intersection of $S$ with $D'$, and let $g_y$ be the geodesic
segment in the support of $\lambda$ centered at $y$ and of length $2\epsilon_0$. 
Then the angle between $g_y$ and $D'$ at $y$ is at least $\pi/4$.
\end{claim}

\begin{proof}
We call $g_x$ the geodesic segment contained in the support of $\lambda$ centered
at $x$ and of length $2\epsilon_0$. $g_x$ is disjoint from $g_y$ on $S$ while 
$x$ is at distance at most $\rho_0$ from $y$, it 
follows that there exists $c>0$ (depending on $\epsilon_0$ and $\rho_0$, and
going to $0$ as $\rho_0\rightarrow 0$ for fixed $\epsilon_0$) such that 
the distance to $g_x$ in $S$ of any point of $g_y([-\epsilon_0+\rho_0,\epsilon_0-\rho_0])$ 
is at most $c\epsilon_0$. 
The same estimate holds in $\Mt$, the universal cover of $M$. If $\rho_0$ is
small enough --- relative to $\epsilon_0$ --- the result follows. 
\end{proof}

\begin{remark} \label{rk:k}
There exists $k_0>0$, depending only on $\rho_0$, such that, if $\Omega$ is a 
convex subset in the disk of radius $\rho_0$ in $H^2$, the total curvature
of the boundary of $\Omega$ is at most $k_0$.
\end{remark}

\begin{proof}
This follows from the Gauss-Bonnet Theorem applied to $\Omega$, with 
$k_0$ equal to $2\pi$ plus the area of the hyperbolic disk of radius $\rho_0$.  
\end{proof}

\begin{proof}[Proof of Lemma \ref{lm:max}]
If $\gamma$ does not intersect the support of $\lambda$ the statement obviously
holds, so we suppose that some point $x\in \gamma$ is in the support of $\lambda$.
Let $D$ be the totally geodesic disk of radius $\epsilon_0/2$ centered at $x$ and
orthogonal, at $x$, to the support of $\lambda$. By construction $D$ is disjoint
from $M_s$. 

Remark \ref{rk:k} shows that the total curvature of the connected component 
$c$ of $D\cap S$ containing $x$ is at most $k_0$. By Claim \ref{cl:angle},
each geodesic in the support of $\lambda$ which intersects $c$ makes with $D$
an angle at least $\pi/4$. It follows that $i(c,\lambda)\leq 2k_0$. It also
follows, since the length of $\gamma$ is less than $\epsilon_0/4$,
that $\gamma$ can be deformed transversally to $\lambda$ to a segment of $c$,
so that $i(\gamma,\lambda)\leq i(c,\lambda)$. Therefore $i(\gamma,\lambda)\leq
2k_0$, and this proves the lemma.
\end{proof}

\begin{proof}[Proof of Proposition \ref{pr:average}]
Notice first that Lemma \ref{lm:max}, although stated only for geodesic
segments $\gamma$ that are at distance at least $\epsilon_0/2$ from the
cone singularities, actually applies without this hypothesis. This is because,
by Lemma \ref{lm:dist2}, the support of $\lambda$ cannot enter the 
$\epsilon_0$-neighborhood of the singular points, so that any part of 
$\gamma$ at distance less than $\epsilon_0$ from the singular set of $S$
has zero intersection with $\lambda$.

Let $n\in \N$ be the unique integer such that $n\epsilon_0/4\leq l_m(\gamma)<
(n+1)\epsilon_0/4$. Then $\gamma$ can be cut into a sequence of segments $\gamma_1,
\cdots, \gamma_n$ of length $\epsilon_0/4$ and one last segment $\gamma_{n+1}$
of length smaller than $\epsilon_0/4$. Lemma \ref{lm:max} can be applied
to each of those segments, it yields that $i(\lambda, \gamma_i)\leq \lambda_1$,
$1\leq i\leq n+1$, so that $$ i(\lambda,\gamma) \leq (n+1)\lambda_1 \leq 
\left(\frac{4l_m(\gamma)}{\epsilon_0}+1\right) \lambda_1~, $$
this proves the proposition.
\end{proof}

\subsection{The grafted metric and the hyperbolic metric at infinity.}

We consider here the relation between the grafted metric $G_\lambda(m)$
and the hyperbolic metric at infinity $t$. 

We first recall the definition of the grafting map 
\begin{eqnarray*}
  G: & \cML\times \cT & \rightarrow \cT \\
  & (l,m) & \mapsto G_l(m)
\end{eqnarray*}
on a closed surface $S$.

The definition of $G_l(m)$ is simpler when $l$ is a weighted multicurve,
that is, its support is a disjoint union of closed curves
$c_1, \cdots, c_N$. The transverse measure is then described by a positive
weight $w_i$ on $c_i$, for $1\leq i\leq N$. Then $G_l(m)$ is obtained
by realizing each curve $c_i$ as a closed geodesics in $(S,m)$, cutting
$S$ open along each $c_i$, and gluing in a flat strip of width $w_i$. 
Thurston showed that this map extends by continuity from weighted multicurves
to measured laminations, see \cite{kulkarni-pinkall}.

We now return to the setting where $m$ and $\lambda$ are the induced metric
and measured bending lamination of the convex core of a quasifuchsian manifold $M$.
The grafted metric $G_\lambda(m)$ is then isometric to the induced metric on the
unit normal bundle of $C(M)$ in the unit tangent bundle of $M$. (The unit normal
bundle of $C(M)$ is the space of unit vectors at points of $\partial C(M)$ which
are the oriented normals of a support plane of $C(M)$.) 

One of the key properties of the grafted metric $G_\lambda(m)$ 
(see e.g. \cite{kulkarni-pinkall}) is that it is in the conformal
class $\tau$ at infinity --- more precisely, there is a natural ``Gauss
map'' defined from the unit normal bundle of $\dr C(M)$ with its ``grafted metric'' to the 
boundary at infinity of $M$, which is conformal. This means that 
$G_\lambda(m)$ is conformal to $t$.
Moreover, since the angles $\theta_i$ are in $(0,\pi)$, the intersection
of the boundary of the convex core with the particles is at non-zero
distance (for $m$) from the support of $\lambda$, so that the cone angles
of the grafted metric at the intersections with the particles of the
boundary of the convex core is equal to the cone angle of the corresponding
singularities. 

The fact that $t$ is conformal to $G_\lambda(m)$ translates as
$$ t = e^{2u} G_\lambda(m)~, $$
where $u:\dr M\rightarrow \R$ is a function.

\begin{lemma} \label{lm:hyperbolique}
The function $u$ is non-positive on $\dr M$.
\end{lemma}

\begin{proof}
Consider two metrics $g$ and $g'$ with $g'=e^{2u}g$, and let $K$
and $K'$ be their curvatures. Then (see e.g. chapter 1 of \cite{Be})
$$ K' = e^{-2u}(\Delta u + K)~. $$
We can apply this formula here with $g=G_\lambda(m)$ and $g'=t$,
so that $K'=-1$ while $K\in [-1,0]$. It takes the form:
$$ \Delta u = -K -e^{2u} = |K| - e^{2u}~, $$
with $K\in [-1,0]$ (this equation is understood in a distributional sense). 

Since the cone angles are the same for $t$ and for $G_\lambda(m)$, 
$u$ is continuous and bounded at the singular points (see \cite{troyanov}).
Let $x_M\in S$ be a point where $u$ achieves its maximum. Suppose first that
$x_M$ is not a singular point, then $u$ is $C^2$ at $x_M$ by elliptic
regularity (see \cite{troyanov}). Moreover $\Delta u\geq 0$ at $x_M$ since
$x_M$ is a maximum of $u$. It follows that
$e^{2u}\leq |K|\leq 1$, so that $u\leq 0$. To complete the proof it is 
sufficient to prove that $u$ cannot achieve a positive maximum at a singular point
of $S$. So we consider a singular point $x_0$ of $S$, and suppose that 
$u>0$ at $x_0$. We will show that $u$ cannot have a maximum at $x_0$.

Let $D$ be the geodesic disk of
radius $r$ centered at $x_0$ in $(S,G_\lambda(m))$. Since $\lambda$ does not
enter a small neighborhood of $x_0$, $D$ is hyperbolic, with only one
cone singularity at $x_0$, if $r$ is small enough. Let $i_0$ be the 
isometric map between $D$, with the metric $G_\lambda(m)$, 
and the hyperbolic disk $H^2_\alpha$ 
with one cone singularity of angle $\alpha$, where $\alpha$ is the
cone angle of $S$ at $x_0$.
Let $i_1:D\rightarrow H^2_\alpha$ be the isometric 
embedding of $(D,t)$ in $H^2_\alpha$. Call $v_0$ the vertex of $H^2_\alpha$,
i.e., its singular point. Since $u>0$ at $x_0$, if $r$ is small enough,
then
$$ \forall x\in D\setminus \{ x_0\}, d(i_1(x),v_0)>d(i_0(x),v_0)~. $$

There is a natural complex map $\phi:H^2_\alpha\rightarrow H^2$, given
in holomorphic coordinates centered at the singular point by 
$z\rightarrow z^{2\pi/\alpha}$. It is conformal and multiplies the metric
by a factor $(2\pi/\alpha)^2 d(x,v_0)^{2(2\pi/\alpha-1)}$. 
Consider the composition
$$ \Phi := \phi \circ i_1\circ i_0^{-1} \circ \phi^{-1}:(\phi\circ i_0)(D)
\rightarrow (\phi\circ i_1)^{-1}(D)~. $$
It is a conformal map, with conformal factor equal to 
$$ (2\pi/\alpha)^2 d(i_1(x),v_0)^{2(2\pi/\alpha-1)} e^{2v} (2\pi/\alpha)^{-2} 
d(i_0(x),v_0)^{-2(2\pi/\alpha-1)}~, $$
with $v=u\circ (\phi\circ i_0)^{-1}$. 
This can be written as
$$ \left(\frac{d(i_1(x),v_0)}{d(i_0(x),v_0)}\right)^{2(2\pi/\alpha-1)}e^{2v}~, $$
and is bigger than $1$ since $u>0$ and $ d(i_1(x),v_0)>d(i_0(x),v_0)$.

Since $\Phi$ is a conformal map 
between two hyperbolic domains, its conformal factor cannot have a 
local maximum bigger than $1$ at an interior point 
by the argument given at the beginning
of this proof. Therefore, $u$ cannot have a positive maximum at $x_0$.
\end{proof}

The following notion will be useful in this section and the next. 

\begin{df}
A {\it c-curve} on $\dr M$ is either a closed curve or a segment with
endpoints at cone singularities, which does not contain any singular
point (except at its endpoints if it's not a closed curve).
\end{df}

We will sometimes implicitly consider c-curves up to homotopy in
the complement of the singular points in $\dr M$. Each homotopy
class (with fixed endpoints) contains a unique geodesic for any non-positively
curved metric on $\dr M$ (in particular for $m, t$ and $G_\lambda(m)$).
Given a c-curve $\gamma$, we will denote by $L_m(\gamma)$ (resp.
$L_t(\gamma)$, $L_{G_\lambda(m)}(\gamma)$) the length of that geodesic
for the corresponding metric.

\begin{cor} \label{cr:hyperbolique}
Let $\gamma$ be a c-curve in $\dr M$, then $L_t(\gamma) \leq 
L_{G_\lambda(m)}(\gamma)$.
\end{cor}

This follows directly from Lemma \ref{lm:hyperbolique}, since any minimizing
c-curve in $(S, G_\lambda(m))$ has shorter length for $t$, and the
minimizing curve in $(S,t)$ in the same homotopy class is even shorter.
Note also that for any $c$-curve $\gamma$, 
$i(\lambda,\gamma)\leq CL_m(\gamma)$. This follows from Proposition \ref{pr:average},
and by the fact that the lengths of the $c$-curves which are segments between
two singular points of $S$ is bounded from below. 

\subsection{An upper bound on the lengths of the curves at infinity}

The second step in the proof of Proposition \ref{pr:length} is a comparison 
between the lengths of c-curves in the metrics $t$ and $m$. 

\begin{prop} \label{pr:longueur-maj}
There exists a constant $C>0$ (independent of $M$) such that:
\begin{enumerate}
\item for each c-curve $\gamma$ in $\dr M$,
$L_{t}(\gamma) \leq C L_{m}(\gamma)$,
\item for each long tube $T$ in the thin part of $(\dr M, m)$, $T$ might
also be a long tube for $t$, but its length for $t$ is at most its
length for $m$ plus $C$.
\end{enumerate}
\end{prop}

The proof uses some simple statements on the geometry of long
hyperbolic tubes in $(S,t)$. Recall (see \cite{dryden-parlier})
that the Margulis Lemma applies to hyperbolic surfaces with cone
singularities of angle at most $\theta$, when $\theta\in (0,\pi)$: 
there exists a constant $c_M$, depending on $\theta$ only, such 
that the set of points where the injectivity radius is less than
$c_M$ is a disjoint union of cusps, disks centered at a cone
singularity, and tubes with core of length less than $2c_M$.

We consider in this subsection a hyperbolic tube
$T$, which can be described as isometric to the set of points at distance at
most $L$ (for some $L>0$) from the unique simple closed geodesic
in the quotient of the hyperbolic plane $H^2$ by a hyperbolic
translation of length $l$. Moreover $l$ is supposed to be small
and $L$ large, so that the lengths of the boundary components 
of $T$ --- which are both equal to $l\cosh(L)$ --- are equal to $c_M$. 
We call $\sigma$ the core of $T$, in other terms the unique
simple closed geodesic contained in $T$, and we denote by
$\sigma_M$ the $c_M$-neighborhood of $\sigma$ --- the set of points 
at distance at most $c_M$ from $\sigma$ in $T$.

\begin{lemma} \label{lm:intersection}
There exists a constant $C>0$ such that $i(\lambda, \dr\sigma_M)\leq
Ce^{-L}$:
the intersection of $\lambda$ with the boundary of $\sigma_M$ is at
most $Ce^{-L}$.
\end{lemma}

\begin{proof}
Any maximal embedded geodesic segment in $T$ intersects exactly 
once $\sigma$, but also each of the two connected components of
$\dr \sigma_M$. It follows that the intersection with $\lambda$
of each of the connected components of $\dr\sigma_M$ is equal
to $i(\lambda, \sigma)$. But since the length of $\sigma$
is $l=e^{-L}$, Proposition \ref{pr:average} --- applied to long
segments that wrap many times around $\sigma$ --- shows that
$i(\lambda, \sigma)\leq C e^{-L}$.  
\end{proof}

\begin{lemma} \label{lm:projection}
There exists a constant $C>0$ such that, if $g$ is an embedded
maximal geodesic segment in $T$, then the length of the orthogonal
projection on $\sigma_M$ of $g\cap (T\setminus \sigma_M)$ is at
most $C$.
\end{lemma}

\begin{proof}
If $g\subset \sigma$, then $g\cap (T\setminus \sigma_M)=\emptyset$
and the result applies. We suppose from here on that $g$ is not 
contained in $\sigma$. If $g$ is contained in one connected component of
the complement of $\sigma$ in $T$ then, since $g$ is embedded, its orthogonal
projection on $\sigma$ is injective, so that the length of its orthogonal
projection is bounded by the length of $\sigma$. Otherwise, 
it follows from standard hyperbolic geometry arguments that $g$ 
intersects $\sigma$ exactly once.

Consider the universal cover $\Tt$ of $T$, it is isometric to the set of
points at distance at most $L$ from a geodesic $\sigmat\subset H^2$
which is the lift of $\sigma$. Choose one of the connected components,
say $\gt$, of the lift of $g$ to $\Tt$. It intersects the lift of
$\dr \sigma_M$ with an angle which is bounded from below --- otherwise
$g$ could not intersect $\sigma$. It follows from this, and from 
elementary geometric properties of the hyperbolic plane, that the 
length of the 
orthogonal projection on $\sigmat$ of each of the segment of 
$\gt$ outside the lift of $\sigma_M$ is bounded from above by
a constant.
\end{proof}

\begin{cor} \label{cr:intersection}
There exists a constant $C>0$ such that, whenever $g_0$ is a maximal
geodesic segment in $T$ such that the orthogonal projection of $g_0$
on $\sigma$ is injective, then 
$i(\lambda_{|T\setminus \sigma_M}, g_0)\leq C$.
\end{cor}

\begin{proof}
Let $c$ be a maximal geodesic segment in the intersection with $T$ 
of the support of $\lambda$,
and let $c'$ be one of the connected components of $c\cap T\setminus
\sigma_M$. Since both $c$ and $g_0$ are geodesic segments, the
union of the orthogonal projections on $\sigma$ of the segments of
$g_0$ and of $c$ between two successive intersections between
them covers $\sigma$. 

It follows that the number of intersections between $c'$ and $g_0$ is
at most equal to $(l_{c'} + l_{g_0})/l$, where $l_{c'}$ is the
length of the orthogonal projection of $c'$ on $\sigma$ and 
$l_{g_0}$ is the length of the orthogonal projection of $g_0$
on $\sigma$ (and $l$ is the length of $\sigma$). 

But the hypothesis on $g_0$ shows that $l_{g_0}\leq l$, while 
Lemma \ref{lm:projection} shows that $l_{c'}\leq C$.
So the number of intersections between $c'$ and $g_0$ is at
most $Ce^{L}$, where $C$ is some positive constant.

Since this inequality applies to all geodesic segments in the
support of $\lambda$, we find that
$$ i(\lambda_{T\setminus \sigma_M}, g_0) \leq Ce^L i(\lambda, 
\dr \sigma_M)~, $$ 
and Lemma \ref{lm:intersection} then shows that 
$i(\lambda_{T\setminus \sigma_M}, g_0)$ is bounded by a positive
constant.
\end{proof}

\begin{proof}[Proof of Proposition \ref{pr:longueur-maj}]
Let $\gamma$ be a c-curve in $\dr M$. It follows from Proposition
\ref{pr:average} that
$$ L_{G_\lambda(m)}(\gamma) \leq L_m(\gamma) + i(\lambda, \gamma) \leq 
C(L_m(\gamma)+1)~, $$
where here again $C$ is a constant depending only on the topology
of $M$. Moreover Corollary \ref{cr:hyperbolique} indicates that
$$ L_t(\gamma)\leq L_{G_\lambda(m)}(\gamma)~, $$
and point (1.) follows.

For point (2.) consider a closed geodesic $\gamma$ contained in the union 
of $T$ and of the thick part of $\dr M$, such that
\begin{itemize}
\item the intersection of $\gamma$ with $T$ has two connected components
$\gamma_1$ and $\gamma_2$,
\item the intersection of $\gamma$ with the thick part of 
$\dr M$ (for $m$) has two connected components $\gamma'_1$ 
and $\gamma'_2$, and each has length bounded by $C$.
\end{itemize}
If $T$ separates the boundary component of $\dr M$ 
containing it, $\gamma$ has to go through $T$ twice, otherwise it is
not necessary but it is still possible to choose $\gamma$ with this
property, and both cases can then be treated in a uniform manner.

Once such a curve $\gamma$ has been found, it is possible to change it
by Dehn twists
so that, in addition to the conditions above, the segments $\gamma_1$
and $\gamma_2$ ``wrap'' at most once around $T$, i.e., their orthogonal
projection to the core $\sigma$ of $T$ is injective. This is achieved
by ``untwisting'' $\gamma$ as much as is necessary. 

Denote as above by $\sigma$ the core of $T$, and by $\sigma_M$
the set of points at distance at most $c_M$ from $\sigma$. Since 
$\gamma$ wraps at most once around $T$, the length of the
intersection of $\gamma_1$ and $\gamma_2$ with $\sigma_M$ is at
most $3c_M$. It then follows from Proposition \ref{pr:average} that 
$$ i(\lambda_{|\sigma_M}, \gamma)\leq C~, $$
where $C$ is some positive constant. By the same proposition, the
intersection of $\gamma$ with the restriction of $\lambda$
to the thick part of $(\dr M, m)$ is at most $C$. 
But Corollary \ref{cr:intersection} shows that
$$ i(\lambda_{|T\setminus \sigma_M}, \gamma)\leq C~. $$

Putting together those estimates we obtain that $i(\lambda, 
\gamma)\leq C$, where $C$ is yet another positive constant. 
The definition of the grafted metric then proves that 
$$ L_{G_\lambda(m)}(\gamma) \leq L_{m}(\gamma) + C~. $$
Finally Lemma \ref{lm:hyperbolique} indicates that the
length of $\gamma$ for $t$ is less than that for 
$G_\lambda(m)$. The result follows.
\end{proof}

\subsection{A bound on the length of the earthquake lamination}
\label{sc:length}

We now switch from 3-dimensional to 2-dimensional geometry
to show that an upper bound on the length of curves in 
$(\dr M,t)$ --- relative to the length of the same curves
in $(\dr M, m)$, as stated in Proposition \ref{pr:longueur-maj} 
--- implies a lower bound on the same lengths. Proposition
\ref{pr:length} will follow. 

We consider a closed
surface $\Sigma$, with some marked points $x_1, \cdots,
x_{n_0}$, and an angle $\theta_i\in (0,\pi)$ attached to $x_i$.

\begin{prop} \label{pr:max-length}
For each $C>0$ there is a constant $C'>0$ as follows.
Let $h,h'\in \cH_{\Sigma, x,\theta}$ be two hyperbolic metrics
such that:
\begin{enumerate}
\item for each c-curve $\gamma$ in $\Sigma$,
$L_{h'}(\gamma) \leq C L_h(\gamma)$,
\item for each long tube $T$ in the thin part of $(\Sigma, h)$, $T$ might
also be a long tube for $h'$, but its length for $h'$ is at most its
length for $h$ plus $C$.
\end{enumerate}
Let $\nu\in \cM_{\Sigma,x}$ be the measured lamination such that 
$h' = E_r(\nu)(h)$. Then the length $L_h(\nu)$ is at most
$C'$.
\end{prop}

The proof of Proposition \ref{pr:max-length}
will use a basic estimate on the variation of the 
length of curves under an earthquake, essentially taken from
\cite{cone}.

\begin{prop} \label{pr:estimee}
Let $m\in \cML_{\Sigma,x}$ be a measured lamination, let $g\in 
\cH_{\Sigma,x,\theta}$ be a hyperbolic metric with cone singularities,
and let $g':=E^r_m(g)$.
Let $\gamma$ be a c-curve. Then
$$ |L_g(\gamma) - L_{g'}(\gamma)|\leq 
i(m,\gamma) \leq L_g(\gamma) + L_{g'}(\gamma)~. $$
\end{prop}

\begin{proof}
The upper bound on $i(m,\gamma)$ can be found in \cite{cone} 
(Lemma 7.1, p. 76); it is stated there for closed curves, but
the proof extends directly to segments between two singular 
points. 

For the lower bound on $i(m,\gamma)$, suppose first that the
support of $m$ is a disjoint union of simple closed curves.
Consider the geodesic (for $g$) $\gamma_0$ homotopic to 
$\gamma$ in $\Sigma_x$, and let $\gamma_1$ be its image by the 
earthquake $E^r_m$, along with the union of the segments in 
the support of $m$ between two points corresponding --- after the 
earthquake --- to one intersection of $m$ with $\gamma_0$.
$\gamma_1$ is homotopic to $\gamma_0$ in $\Sigma_x$.
Clearly $L_{g'}(\gamma_1)=L_g(\gamma) + i(m,\gamma)$, 
while $L_{g'}(\gamma)\leq L_{g'}(\gamma_1)$. It follows that
$L_{g'}(\gamma) \leq L_g(\gamma) + i(m,\gamma)$. The same 
inequality also holds with $g$ and $g'$ exchanged, and
the lower bound on $i(m,\gamma)$ follows. The result when 
$m$ is a general lamination --- not rational --- holds by
density of the rational laminations in $\cML_{\Sigma,x}$.
\end{proof}

We now return to the notations used in Proposition \ref{pr:max-length}.
Note that the support of $\nu$ is a geodesic lamination in
$(\Sigma, g)$. It is therefore possible to consider the 
intersection of $\nu$ with the thin (resp. thick) part of $\Sigma$ for 
$g$, which we call $\nu_t$ (resp. $\nu_T$). The same 
decomposition can be done for $g'$, leading to $\nu'_t$
and $\nu'_T$. 

We first state a basic property of hyperbolic surfaces, which
is necessary below.

\begin{lemma}
There exist $r>0$, $C>0$ and $\theta_0\in (0,\pi)$, depending
only on the supremum $\theta_M$ of the $\theta_i$ and on the genus of $\Sigma$, 
such that, for any $x\in \Sigma_T$ and any geodesic segment $\gamma_0$
of length $2r$ centered at $x$, there exists a closed geodesic in 
$\Sigma$ of length at most $C$ intersecting $\gamma_0$ with angle
at least $\theta_0$.
\end{lemma}

\begin{proof}
Note that any maximal segment in the thick part of a topologically 
non-trivial hyperbolic surface (with cone singularities of angle less
than $\pi$) intersects some closed geodesic, of length bounded by
a constant $C$. 

The statement therefore
follows from a straightforward compactness argument. Indeed, if the
constant $\theta<\theta_M$ did not exist, there would be a sequence of
thick hyperbolic surfaces with boundary $\Sigma_{T,n}$ 
(with cone singularities of angles less
than $\pi$), for which the optimal value of $r$ would go to infinity, or the optimal
value of $\theta$ would go to $\pi$, as $n\rightarrow \infty$. This sequence
could be taken of fixed topology, and the diameter of those surfaces would 
then be bounded, so that $r$ would necessarily be bounded.
 
We could then choose a converging subsequence, and obtain a thick hyperbolic
surface (with cone singularities of angle less than $\theta_M$) for which
some maximal geodesic segment intersects no closed geodesic of length 
less than $C$ transversally, a contradiction.
\end{proof}

\begin{lemma} \label{lm:long-epais}
There exists a constant $C$ (depending only on the genus of 
$\Sigma$) so that the length of $\nu_T$ is at most $C$.
\end{lemma}

\begin{proof}
Let $r_0>0$ be smaller than the injectivity radius of $(\Sigma,
h)$ at each point of $\Sigma_T$. There exists another number 
$r_1\in (0,r_0)$ with the following property: if $\gamma_0$ and 
$\gamma_1$ are two disjoints geodesics in $H^2$ and $x\in \gamma_0$
is at distance at most $r_1$ from $\gamma_1$, then any geodesic
intersecting $\gamma_0$ at distance less than $r_1$ from $x$
and making an angle bigger than $\theta_0$ with 
$\gamma_0$ intersects $\gamma_1$ at distance at most $r_0$ 
from $x$.

Choose a large constant $C_1>0$.
If the length of $\nu_T$ were bigger than some large constant,
the sum of the weights of the segments of the support of 
$\nu$ intersecting some geodesic disk of radius $r_1$ and center
$x\in \mbox{supp}(\nu)\cap \Sigma_T$ would be bigger than $C_1$. 
Applying the previous lemma, with $\gamma_0$ equal to a
segment containing $x$ in the support of $\nu$, would yield 
a closed curve $c$ in $\Sigma_T$, of bounded length, such
that $i(c,\nu)$ is arbitrarily large. 

Proposition \ref{pr:estimee} would then show that the length of 
$c$ for $h'$ is much larger than the length of $c$ for $h$, 
contradicting point (1) in the hypothesis of Proposition
\ref{pr:max-length}.
\end{proof}

\begin{lemma} \label{lm:long1}
There exists a constant $C>0$ as follows. Let $\gamma\subset
\Sigma_T$ be a geodesic segment of length at most $c_M$. Then 
$$ i(\nu,\gamma)\leq CL(\nu_T)~. $$
\end{lemma}

\begin{proof}
We will consider the case when $\nu$ is rational, the general case
follows by density of the rational measured laminations.

Let $r$ be the injectivity radius of $\Sigma_T$. Let $\nu_\gamma$
be the union of the intersections with $\Sigma_T$ of all 
geodesic segments centered at a point $x\in\gamma$, of length 
$2r$, in the support of $\nu$. Each of those segments has 
length at least $r$, since at least one side of $x$
is contained in $\Sigma_T$. 

By definition of $r$,
those segments intersect $\gamma$ exactly once. Moreover, the
length of $\nu_T$ is larger than the sum over the segments of their
length times their weight (this sum is finite since $\nu$ is
rational). But this sum is at least $ri(\nu, \gamma)$, so that
$L(\nu_T)\geq ri(\nu, \gamma)$. This proves the lemma.
\end{proof}

\begin{lemma} \label{lm:poids-ame}
There exists a constant $C>0$ such that,
if $T$ is a tube of length $2L$ in $(\Sigma_t, h)$, with 
core $\sigma$, then $i(\nu,\sigma)\leq Ce^{-L}$.
\end{lemma}

\begin{proof}
$L_h(\sigma) = c_1e^{-L}$, where $c_1$ is some constant.  
Point (1) in the hypothesis of Proposition \ref{pr:max-length} shows that 
the length of $\sigma$ for $h'$ is at most $c_2e^{-L}$, where 
$c_2$ is another positive constant. But Proposition \ref{pr:estimee}
then yields the result.
\end{proof}

Recall that %% In the next lemma, 
we call $\sigma_M$ the set of 
points at distance at most $c_M$ from $\sigma$.

\begin{lemma} \label{lm:poids-vois}
There exists  a constant $C>0$ such that,
if $T$ is a tube of length $2L$ in $(\Sigma_t, h)$, with 
core $\sigma$, then the length for $h$ of the restriction of 
$\nu$ to $\sigma_M$ is at most $C$.
\end{lemma}

Note that there is no reason to believe that this statement %% proposition
is optimal; indeed, it appears quite reasonable to think that
the bound could be improved to $Ce^{-L}$. The bound given
here, however, is both sufficient for our needs and easier to
obtain.

\begin{proof}
We know by Lemma \ref{lm:poids-ame} that $i(\nu,\sigma)\leq
C_1e^{-L}$, where $C_1>0$ is some constant. It follows that, if 
$L_h(\nu_{|\sigma_M})$ is larger than some constant $C_2$, then 
each leaf of $\nu$ intersects $\sigma$ with a very small angle,
and for any geodesic segment $\gamma$ going through $T$ and
intersecting $\sigma$ with angle bigger than $\pi/4$, $i(\gamma,
\nu_{|\sigma_M})\geq C_2e^L$.

Let $\gamma_1$ be a closed geodesic in $(\Sigma,h)$ which has two 
segments in $\Sigma_T$ and two segments going through $T$. Furthermore we choose $\gamma_1$ with the smallest length. 
Since $\Sigma_T$ has bounded diameter, there is $C>0$ such that $\ell_h(\gamma_1)\leq C+2T$. Since $\gamma_1$ has minimal length, it intersects
$\sigma$ with angle at most $\pi/4$ at each of the two intersections.
Then $i(\nu, \gamma_1)\geq 2C_2e^L$, so that, by Proposition 
\ref{pr:estimee}, $L_{h'}(\gamma_1)$ is much larger than 
$L_h(\gamma_1)$. This contradicts the hypothesis of Proposition
\ref{pr:max-length}.
\end{proof}

\begin{proof}[Proof of Proposition \ref{pr:max-length}]
According to Lemma \ref{lm:long-epais}, the length of
the restriction of $\nu$ to the thick part of $\Sigma$ is 
bounded by a constant (depending only on the genus of $\Sigma$).
$\Sigma$ is the union of $\Sigma_T$, a finite set of neighborhoods of
cusps and cone singularities (which can be disregarded because geodesic
laminations do not enter them), and a finite set of 
long thin tubes, the number of those tubes being at most 
$3g-3$, where $g$  is the genus of $\Sigma$. 
Let $T$ be one of those tubes, and let $\sigma$ be its
core. Then $L(\nu_{|\sigma_M})$ is bounded by a constant
by Lemma \ref{lm:poids-vois}. Moreover, the length
of each maximal segment of the support of $\nu$ in $T\setminus
\sigma_M$ is at most $2e^L$, and each is contained in a 
maximal segment in $T$ (contained in the support of $\nu$)
which intersects $\sigma$ once. Since $i(\nu,\sigma)\leq
Ce^L$, the length of the restriction of $\nu$ to $T\setminus
\sigma_M$ is at most $4C$. Summing all contributions to the
length of $\nu$ yields the desired result.
\end{proof}

\begin{proof}[Proof of Proposition \ref{pr:length}]
The statement clearly follows from Proposition \ref{pr:longueur-maj} 
and from Proposition \ref{pr:max-length}.
\end{proof}

%% file: conebend6.tex
\section{The conformal structure at infinity}
\label{se:induced}

This section contains the proof of Theorems \ref{tm:bers} and \ref{tm:metriques}, 
mostly as a consequence of Lemma \ref{pr:length}. 

\subsection{A topological lemma}

We first state a simple topological lemma, necessary below to apply Proposition
\ref{pr:length} as directly as possible. We fix a closed surface $S$ of genus 
at least $2$, a $n$-tuple
of points $x=(x_1,\cdots, x_{n_0})$ and a $n_0$-tuple of angles $\theta=(\theta_1,\cdots,
\theta_{n_0})\in (0,\pi)^n$.

\begin{lemma} \label{lm:topo}
Let $c>0$, and let $K\subset \cH_{S,x,\theta}$ be a compact subset. The set of all elements
of $\cH_{S,x,\theta}$ obtained by a right earthquake along a measured lamination of length
at most $c$ on an element of $K$ is relatively compact.
\end{lemma}

\begin{proof}
Let $m\in K$. The set of measured laminations $l\in \cML$ of length less than $c$
for $m$ is compact in $\cML_{S,x}$. Since the earthquake map is continuous relative to
the measured lamination factor, the set 
$$ E_r(\{ l\in \cML_{S,x}~|~L_m(l)\leq c\}\times \{m\}) $$ 
is compact in $\cH_{S,x,\theta}$. 

Again because the earthquake map $E_r$ is continuous, it follows that
there is a neighborhood $U_m$ of $m$ in $\cH_{S,x,\theta}$ such that the image by $E_r$ of
$$ \{ (l,m') \in \cML_{S,x}\times U_m ~|~ L_{m'}(l)\leq C\} $$
is relatively compact. 

Since $K$ is compact, it is covered by finitely many such neighborhoods 
$U_{m_i}$, for $m_i$ in $K$. The result follows.
\end{proof}

\subsection{Compactness relative to the conformal structure at infinity} \label{ssc:62}

The previous considerations lead to a simple proof of 
Proposition \ref{pr:compact-conforme} from Proposition \ref{pr:length} and 
Lemma \ref{lm:conv}.

Consider a sequence $(g_n)_{n\in \N}$ of quasifuchsian metrics with particles,
as in Proposition \ref{pr:compact-conforme}. Let $(m_n)_{n\in \N}$ be the
sequence of induced metrics on the boundary of the convex core, and let
$t_n$ be the sequence of hyperbolic metrics in the conformal class at 
infinity $\tau_n$. Since $(\tau_n)_{n\in \N}$ converges by the hypothesis
of Proposition \ref{pr:compact-conforme}, $(t_n)_{n\in \N}$ converges
to a limit $t$, so it remains in a compact subset of $\cH_{S,x,\theta}\times 
\cH_{S,x,\theta}$. 

But $m_n$ is obtained from $t_n$ by an earthquake along a measured lamination
$\nu_n$ which, by Proposition \ref{pr:length}, has bounded length. Lemma 
\ref{lm:topo} therefore shows that $(m_n)_{n\in \N}$ remains in a compact
subset of $\cH_{S,x,\theta}\times \cH_{S,x,\theta}$. We can therefore extract
a sub-sequence so that $(m_n)_{n\in \N}$ converges.

Lemma \ref{lm:conv} then shows that $(g_n)_{n\in \N}$ has a subsequence
which converges to a quasifuchsian metric with particles. This proves the
proposition.

\subsection{Proof of Theorem \ref{tm:bers}}

We are now ready to prove Theorem \ref{tm:bers}. It is helpful to
introduce additional notations: 
\begin{itemize}
\item $\cM_{S,x}:=\cup_{\theta \in (0,\pi)^N} \cM_{S,x,\theta}$ is the 
space of quasifuchsian metrics with particles on $S\times \R$ with
a fixed number of particles but with varying angles,
\item $\cH_{S,x}:=\cup_{\theta \in (0,\pi)^N} \cH_{S,x,\theta}$ is the 
space of hyperbolic metrics on $S$ with a fixed number of cone
singularities but with varying angles,
\item $\Delta_{S,x}:=\cup_{\theta \in (0,\pi)^N} \cH_{S,x,\theta}\times
\cH_{S,x,\theta}$ is a kind of diagonal with respect to the angle
variable in $\cH_{S,x}\times \cH_{S,x}$.
\end{itemize}
Note that, by a result of Troyanov \cite{troyanov} already mentioned 
above, $\cH_{S,x}$ can be naturally identified with $\cT_{S,x}\times
(0,\pi)^N$. The notation is nonetheless useful in the argument 
presented here.

Consider the natural map:
$$ \Phi_{S,x}:\cM_{S,x}\rightarrow \cH_{S,x}\times 
\cH_{S,x} $$
sending a hyperbolic metric with particles on $S\times \R$ (with
cone angles given by the $\theta_i$) to the conformal structures
at $\pm\infty$. It follows from the definition that the image
of $\Phi_{S,x}$ is contained in $\Delta_{S,x}$.

Let $\Phi_{S,x,\theta}$ be the restriction of $\Phi_{S,x}$ to 
$\cM_{S,x,\theta}$, for a fixed $\theta\in (0,\pi)^N$. 
The main result of \cite{qfmp} is that --- in a slightly more
general context, allowing for more topology and for singularities along a graph --- 
$\Phi_{S,x,\theta}$ is a local homeomorphism from $\cM_{S,x,\theta}$ to $\cH_{S,x,\theta}\times 
\cH_{S,x,\theta}$. It follows that $\Phi_{S,x}$ is a local homeomorphism
from $\cM_{S,x}$ to $\Delta_{S,x}$. 
Moreover,  $\Phi_{S,x}$ is proper by Proposition 
\ref{pr:compact-conforme}, so that it is a covering of 
$\cT_{S,x,\theta}\times \cT_{S,x,\theta}$.

To prove that $\Phi_{S,x}$ is a (global) homeomorphism
we need to show that some elements of the target space have
exactly one inverse image. Suppose that for all $i\in \{1,\cdots,
N\}$, $\theta_i=2\pi/k_i$, where $k_i\in \N, k_i\geq 2$. 
Let $\tau_+,\tau_-\in \cH_{S,x,\theta}$. There exists a finite
covering $\pi:\Sb\rightarrow S$, with ramification of order
$k_i$ at the $x_i$, such that $\tau_+$ (resp. $\tau_-$) lifts
to a non-singular hyperbolic metric $\taub_+$ (resp.
$\taub_-$). By the Bers double uniformization theorem \cite{bers}
$\taub_+$ and $\taub_-$ are in the conformal class at infinity
of a unique quasifuchsian hyperbolic metric, say $\gb$, on 
$\Sb\times \R$. Since it is unique, $\gb$ is invariant under
the deck transformations of $\pi$, so that $\gb$ is the 
pull-back to $\Sb\times \R$ of a hyperbolic metric $g$ on 
$S\times \R$, with cone singularities of angle $\theta_i$ along
$\{x_i\}\times \R$, $1\leq i\leq N$. This construction also
shows that $g$ is unique --- since any other hyperbolic metric 
with particles of the given angles would lift to a non-singular
quasifuchsian metric on $\Sb\times\R$, which would have to be $\gb$.
This shows that $(\tau_+,\tau_-)$ has a unique inverse image by
$\Phi_{S,x}$, so that $\Phi_{S,x}$ is a homeomorphism from 
$\cM_{S,x}$ to $\Delta_{S,x}$.

\subsection{Proof of Theorem \ref{tm:metriques}}

We need another natural map.

\begin{df}
Let $\Psi_{S,x,\theta}:\cH_{S,x,\theta}\times \cH_{S,x,\theta}
\rightarrow \cH_{S,x,\theta}\times \cH_{S,x,\theta}$ 
be defined as follows. Given 
$(t_+,t_-)\in \cH_{S,x,\theta}\times \cH_{S,x,\theta}$
and $\theta=(\theta_1,\cdots,\theta_{n_0})\in (0,\pi)^n$, 
there is by Theorem \ref{tm:bers} a unique quasifuchsian 
metric with particles $g\in \cM_{S,x,\theta}$ such that 
$\Phi_{S,x,\theta}(g)=([t_+],[t_-])$. Then
$$ \Psi_{S,x,\theta}(t_+,t_-) = (m_+, m_-)~, $$
where $m_+$ and $m_-$ are the conformal classes of
the induced metrics on the 
upper and lower boundary components of the convex core of
$(S\times \R,g)$, respectively.
\end{df}

According to Proposition \ref{pr:length} and Lemma \ref{lm:topo}, if $(m_+,m_-)=
\Psi_{S,x,\theta}(t_+,t_-)$, then $m_+$ is
$C_q$-quasi-conformal to $t_+$, and $m_-$ is
$C_q$-quasi-conformal to $t_-$. This shows that 
$\Psi_{S,x,\theta}$ is proper and extends continuously to a
map which is the identity on the boundary at infinity of
$\cH_{S,x,\theta}\times \cH_{S,x,\theta}$, so that it is onto.
This proves Theorem \ref{tm:metriques}.

%% file: conebend7.tex
\section{Some questions and remarks}
\label{se:questions}

\subsection{Some questions.}

The quasifuchsian cone-manifolds described above are direct
extensions of the ``usual'' quasifuchsian hyperbolic manifolds which have
received much attention over the last couple of decades. It is quite natural
to wonder whether some properties which have
been conjectured in the smooth case can be extended to the singular setting.

\begin{question}
Does uniqueness hold in Theorem \ref{tm:metriques}?
\end{question}

Another natural question, which is ``dual'' to the previous one in a precise
sense (see e.g. \cite{vbc}) concerns the measured bending lamination on the
boundary of the convex core. 

\begin{question}
Does uniqueness  hold in Theorem \ref{tm:general}?
\end{question}

The same questions can be asked for submanifolds of quasifuchsian
cone-manifolds which are convex but have a smooth boundary, which is
orthogonal to the singular locus. In the smooth case it is known \cite{hmcb}
that one can prescribe the induced metric on the boundary, as well as the its
third fundamental form (the smooth analog of the measured bending lamination)
and that each is obtained uniquely, it would be interesting to know whether
the same is true for quasi-fuchsian cone-manifold. The methods of \cite{hmcb}
do not appear to extend directly to the singular case.

%% ajoute paragraphe suivant
Since the Alhfors-Bers theorem extends as Theorem \ref{tm:bers} to quasifuchsian 
manifolds with particles, it is quite natural to ask whether the Ending Lamination
Conjecture (see \cite{brock-masur-minsky,brock-canary-minsky}) can also be 
extended to hyperbolic manifolds with particles. A natural starting point would
be to consider manifolds homeomorphic to $S\times \R$, where $S$ is a closed
surface of genus at least $2$.

Note also that those questions are not necessarily restricted to
quasifuchsian cone-manifolds, and could also be asked for ``convex co-compact
cone-manifolds'', if that term is understood in a proper way. 

\subsection{AdS manifolds with particles.}

Mess \cite{mess,mess-notes} discovered a remarkable analogy between quasifuchsian
hyperbolic 3-manifolds and globally hyperbolic maximal compact (GHMC) AdS manifolds.
In particular he proved an analog of the Bers double uniformization theorem form
GHMC AdS manifolds: on a manifold homeomorphic to $S\times \R$, where $S$ is a 
closed surface of genus at least $2$, the space of GHMC AdS manifolds is
parametrized by the product of two copies of the Teichm\"uller space of $S$, through
the ``left'' and ``right'' parts of the holonomy representation. 

GHMC AdS manifolds also have a convex core, whose boundary has a hyperbolic
induced metric, as in the quasifuchsian setting, and is pleated along a measured
lamination. The analog of Theorem \ref{tm:general} holds in that context \cite{earthquakes}: any 
two measured laminations that fill $S$ can be obtained as the bending lamination
of the boundary of the convex core. But the uniqueness remains elusive, as in 
the quasifuchsian setting. Moreover, the analog of Theorem \ref{tm:metriques}
is also conjectured to be true but no result is known.

It is also possible to consider GHMC AdS manifolds with ``particles'', i.e., cone
singularities along time-like geodesics, for which the angle is less than $\pi$.
The analog of the Bers double uniformization theorem (more directly, the 
analog of Theorem \ref{tm:bers}) holds in this AdS setting \cite{cone}.
The analog of Theorem \ref{tm:general} is also true in that setting \cite{earthquakes}.
However no analog of Theorem \ref{tm:metriques} is known.

%% paragraphe ajoute
Still in the AdS setting, new phenomena can occur when the singularity is along a graph
(so that the particles are allowed to interact), see \cite{colI,colII}. One can associate to
a GHMC AdS manifold with a graph of interacting particles a sequence of pairs of points
in the Teichm\"uller space of the underlying surface, with each pair corresponding to a
``slice'' where no interaction occurs. It would be interesting to know whether any
analog of this description holds for quasifuchsian hyperbolic manifolds with cone
singularities along a graph (perhaps with some conditions on the cone angles, for instance
cone angles less tha $2\pi$).

\subsection{The renormalized volume.}

Theorem \ref{tm:bers} for quasifuchsian manifolds with particles has
possible applications to the Teichm\"uller theory of hyperbolic surfaces
with cone singularities (of fixed angles) on a surface. Indeed it was
remarked in \cite{volume} that the definition of the renormalized volume of 
a quasifuchsian 3-manifolds extends to manifolds with particles.
Knowing Theorem \ref{tm:bers}, it is possible to remark that the key property
of the renormalized volume --- to be a K\"ahler potential for the Weil-Petersson
metric on Teichm\"uller space --- extends to the natural Weil-Petersson metric
on the Teichm\"uller space of hyperbolic metrics with cone singularities
(of prescribed angle less than $\pi$) on a surface; the proof from \cite{volume,review}
directly extends to this setting. 

One direct consequence is that this Weil-Petersson metric is K\"ahler, as
was discovered by Schumacher and Trapani \cite{schumacher-trapani} by other means.
This metric, however, seems to depend on the choice of the cone angles. 

Another possible application is to some properties of the grafting map considered
on hyperbolic surfaces with cone singularities of angle less than $\pi$, as 
considered in \cite{cp,review}. This is however less directly related to what we
are doing here, since it only uses the geometry of 3-dimensional hyperbolic
ends --- rather than quasifuchsian metrics --- with particles.

%% file: conebend8.tex
\section{Quasiconformal estimates}
\label{se:quasiconf}

This appendix contains the proof of Proposition \ref{pr:qconf}.  
The first step is a simple extension to hyperbolic 
surfaces with cone singularities of some classical tools concerning
pants decompositions.

\subsection{Pants decompositions}

The content of this subsection is probably well known, see e.g. 
\cite{dryden-parlier} for closely related considerations. We include
this material for completeness.

Let $S$ be a closed surface, and let $h$ be a hyperbolic metric
on $S$ with cone singularities at some points $x_1,\cdots, x_{n_0}$,
with cone angles $\theta_1,\cdots,\theta_{n_0}\in (0,\pi)$. If $h$
had cusps --- or geodesic boundary components --- at the $x_i$ rather
than cone singularities, it would be quite natural to consider 
pants decompositions of $(S,h)$. With cone singularities of angles
less than $\pi$, it remains possible.

\begin{df}
A {\it singular pair of pants} is a hyperbolic surface with geodesic boundary,
possibly containing cone singularities of angle less than $\pi$, 
which is either: 
\begin{itemize}
\item a hyperbolic pair of pants (with geodesic boundary) containing 
no cone singularity,
\item a hyperbolic annulus with geodesic boundary containing 
exactly one cone singularity,
\item a hyperbolic disk with geodesic boundary containing exactly
two cone singularities. 
\end{itemize}
\end{df}

Given a singular hyperbolic pair of pants, its three geodesic boundary
components or cone singularities will be called its {\it legs}.
We hope that the reader will excuse us for this weird and 
perhaps confusing terminology. 

\begin{df}
A {\bf pants decomposition} of $S$ is a decomposition 
$S=S_1\cup \cdots \cup S_n$ of $S$ as the union of finitely many
subsurfaces with disjoint interior, each of which is a singular
pair of pants.
\end{df}

It is implicit in this definition that the boundary of the $S_i$
contains no cone singularities; the cone singularities are each 
contained in the interior of one of the singular pairs of pants. 

\begin{lemma} \label{lm:pant-existence}
There exists a constant $C_p>0$ such that, for any choice of $S$ 
and $h$, $(S,h)$ has a pants decomposition with all boundary curves
of length less than $C_p$.
\end{lemma}

\begin{proof}[Sketch of the proof]
A standard recursive argument (see e.g. \cite{benedetti-petronio})
reduces the proof to showing that, for any hyperbolic surface with
cone singularities (of angle less than $\pi$) and geodesic boundary,
there is a simple closed geodesic of length at most $C_p$ which is 
not homotopic to a singular point or to a boundary component.  

This in
turn follows from other standard arguments, for instance  based on
comparing the area of the surface (given by a suitable Gauss-Bonnet
formula, see e.g. \cite{troyanov}) to the area of embedded geodesic
disks.
\end{proof}

\begin{df}
Let $P$ be a singular pant. Its {\it leg invariants} are the
length of its geodesic boundary components and the angles at its
cone singularities. 
\end{df}

For instance, the boundary invariants of a (non-singular) hyperbolic
pair of pants are the lengths of its boundary components. 

\begin{lemma} \label{lm:pant-unique}
Each hyperbolic pair of pants is uniquely determined, up to isometry, 
by its leg invariants and by the type of its ``legs'' --- whether they
are boundary components or cone singularities. 
\end{lemma}

The proof follows the classical arguments used for non-singular
hyperbolic pairs of pants, it is based on elementary properties
of some hyperideal hyperbolic triangles stated below in three 
propositions (the first two have probably 
been known since Lobachevsky). 

Recall that a hyperideal triangle can be defined, using the projective
model of the hyperbolic plane, as a triangle which might have its 
vertices either in the hyperbolic plane, on its ideal boundary, or
outside the closure of the hyperbolic plane (considered as the interior of a disk
in the projective plane), but with all edges intersecting the hyperbolic
plane. A vertex is then ideal if it is on the
ideal boundary, and strictly hyperideal if it is outside the
closed disk. 

Recall also that given a point $v_0$ outside the closure of the projective model of
$H^2$ (in the projective plane), there is a unique hyperbolic geodesic, $v_0^*$,
such that any the intersection with the projective model of $H^2$ of any 
projective line containing $v_0$ is orthogonal to $v_0^*$. This geodesic is
called the line {\it dual} to $v_0$.

We introduce here a slightly restricted notion of hyperideal triangle.

\begin{df}
An {\bf extended hyperbolic triangle} is a hyperbolic triangle with 
one or more strictly hyperideal vertices and its other vertices in
the ``interior'' of the hyperbolic plane. 
A {\bf truncated hyperbolic triangle} is the intersection of 
an extended hyperbolic triangle with the hyperbolic half-planes bounded
by the lines dual to its strictly hyperideal vertices (and not
containing the endpoints of the edges going towards those vertices).
\end{df}

For instance, a right-angle hyperbolic hexagon can be considered ---
in two ways --- as a truncated hyperbolic triangle, with three 
strictly hyperideal vertices. Given a hyperbolic triangle, its 
{\it angles} are the hyperbolic angles at the non-hyperideal vertices 
and the length of its intersections with the lines dual to the 
strictly hyperideal vertices. Notes that those lengths can quite
naturally be considered as angles (they are then imaginary numbers)
but it is not necessary to enter such considerations here (see e.g.
\cite{shu,cpt} for more details). 

There is a natural way to define the {\it edge lengths} of an 
extended hyperbolic triangle. The length of the edge joining two
vertices $v$ and $v'$ is: 
\begin{itemize}
\item the hyperbolic distance between $v$ and $v'$, if neither 
$v$ nor $v'$ is strictly hyperideal,
\item the hyperbolic distance between $v$ and the line dual to 
$v'$, when $v'$ is hyperideal but $v$ is not,
\item the distance between the lines dual to $v$ and $v'$, when
both are strictly hyperideal.
\end{itemize}

It is useful to remark that the lengths and angles of an 
extended hyperbolic triangle satisfy a natural extension of
the cosine formula. Moreover it is quite easy to check that an
extended hyperbolic triangle, with vertices of given type, is
uniquely determined by two lengths and one angle, or by 
two angles and one length. 

\begin{lemma}\label{lm:0}
An extended hyperbolic triangle is uniquely determined by the type of its
vertices --- whether they are ``usual'' or strictly hyperideal
vertices --- and its edge lengths.
\end{lemma}

\begin{proof}
This statement is classical for ``usual'' hyperbolic triangles, with 
no strictly hyperideal vertex. It is also well-known for triangles
with three strictly hyperideal vertices, see \cite{benedetti-petronio}.

Consider an extended hyperbolic triangle, with exactly one strictly
hyperideal vertex, say $v_1$, and two usual vertices, $v_2$ and $v_3$.
Let $l_i$ be the length of the edge between $v_j$ and $v_k$, 
for $\{ i,j,k\}=\{ 1,2,3\}$. Consider $l_2, l_3$ as fixed, then 
$l_1$ is easily seen (using for instance the cosine formula for
extended hyperbolic triangles) to be a strictly increasing function
of the angle $\delta$ at $v_1$. This proves the lemma in this case. 

Consider now the situation where $v_1$ is a ``usual'' vertex, while
$v_2$ and $v_3$ are strictly hyperideal. Then given $v_1$, the positions
of the lines dual to $v_2$ and $v_3$ are completely determined by the
angle $\alpha$ at $v_1$. Moreover the distance between those lines, which
by definition is equal to $l_1$, is a strictly increasing function of
$\alpha$. This shows the result in this last case.
\end{proof}

The same arguments can be used to prove the ``dual'' lemma, concerning
the possible angles. Here we need a more precise statement, in addition
to the fact that the angles determine the triangle we need to know that
a large class of triples of angles can actually be realized. 

\begin{lemma}\label{lm:1}
An extended hyperbolic triangle is uniquely determined by the type of its
vertices --- whether they are ``usual'' or strictly hyperideal
vertices --- and its angles. The angles at the ``usual'' vertices
can take any value in $(0,\pi/2)$, while the ``angles'' at the strictly hyperideal
vertices can be any numbers in $(0,\infty)$.
\end{lemma}

\begin{proof}
Again the case where all vertices are ``usual'' is classical, while
the case where all three vertices are strictly hyperideal is well-known. 

Consider a triangle $T$ with one ``usual'' vertex, say $v_1$, and two
strictly hyperideal vertices, $v_2$ and $v_3$. Let $e_1, e_2, e_3$
be the edges opposite to $v_1, v_2, v_3$ respectively.
The triangle $T$ is completely determined
by the length $l_1$ of the edge $e_1$ and 
by the ``angles'' $\alpha_2$ and $\alpha_3$, that is, 
the lengths of the segments of the lines $v_2^*, v_3^*$ dual to $v_2, v_3$
between their intersections with $e_1$ and with $e_3$ and $e_2$, respectively.

Given $\alpha_2$ and $\alpha_3$, the possible values of $l_1$ vary between 
a minimal value $l_{1,m}$ and a maximal value $l_{1,M}$. Suppose for instance
that $\alpha_2\geq \alpha_3$, then $l_{1,m}$ corresponds to the situation where
$e_3$ is reduced to a point. Then, after truncation, $T$ corresponds to a 
quadrilateral with $3$ right angles. The last angle, between $v_2^*$ 
and $e_2$, has to be less than $\pi/2$ by the Gauss-Bonnet theorem. This
means that for $l_1$ slightly larger than $l_{1,m}$, $\alpha_3>\pi/2$.
On the other hand, $\alpha_3\rightarrow 0$ as $l_1\rightarrow l_{1,M}$, and
$\alpha_3$ is a decreasing function of $l_1\in (l_{1,m},l_{1,M})$. This proves
the lemma for triangles with two strictly hyperideal vertices. 

Consider now a triangle $T'$ with one strictly hyperideal vertex, say $v_1$, and
two ``usual'' vertices, $v_2$ and $v_3$. Consider $\alpha_2,\alpha_3\in (0,\pi/2)$
as fixed, $T'$ is then entirely determined by $l_1$. $l_1$ can vary between a
minimal value $l_{1,m}>0$ and $\infty$, where $l_{1,m}$ corresponds to the case
where $v_1$ is an ideal vertex. The angle $\alpha_1$ then varies between $0$
and $\infty$, and is a strictly increasing function of $l_1$. The result follows.
\end{proof}

\begin{lemma}\label{lm:2}
Each singular pair of pants has a unique decomposition as the 
union of two copies of a truncated hyperbolic triangle (glued
along their common boundary).  
\end{lemma}

\begin{proof}
Let $v_1, v_2, v_3$ be the three legs --- which could be either singular points
or boundary components. There is a unique homotopy class of embedded 
segment joining $v_i$ to $v_j$, for $i\neq j$. Those three segments can be
uniquely realized as minimizing geodesics, which are then orthogonal to the
boundary components. Cutting the pair of pants along those three geodesic
segments yields two extended hyperbolic triangles, glued along their edges. 
Those two extended triangles have the same edge lengths, so that they are
isometric by  Lemma \ref{lm:0}.
\end{proof}

\begin{proof}[Proof of Lemma \ref{lm:pant-unique}]
By Lemma \ref{lm:1}, the two extended triangles glued to obtain a 
hyperbolic pair of pants are uniquely determined by their angles, which can 
take any value as long as the angles at the ``usual'' vertices are less than
$\pi/2$. This shows that hyperbolic pairs of pants are uniquely determined
by their leg invariants, and any values are possible as long as the 
angles at the singular points are less than $\pi$.
\end{proof}

We now turn to the parameterization of hyperbolic metrics with 
cone singularities by Fenchel-Nielsen type coordinates. We first
state a lemma on the existence and uniqueness of a pants decomposition
from topological data, leaving the proof to the reader since it is
the same as in the non-singular case. 

\begin{lemma}\label{lm:pant-topo}
A pants decomposition is uniquely determined by the choice of 
the boundary curves $\gamma_1, \cdots, \gamma_N$, 
considered as simple closed curves in 
$S\setminus \{ x_1,\cdots, x_{n_0}\}$, under the hypothesis that:
\begin{itemize}
\item the $\gamma_i$ can be realized as pairwise disjoint curves,
\item each connected component of their complement is either a 
pair of pants containing none of the $x_i$, or a cylinder containing
exactly one of the $x_i$, or a disk containing exactly two of the
$x_i$.
\end{itemize}
\end{lemma}

Finally we state the main consequence, on the parameterization of
hyperbolic metrics with cone singularities of fixed angle by 
Fenchel-Nielse coordinates, again leaving the proof to the reader.
Note that the Dehn twist parameters are defined only in a relative
way, however this is exactly the same as in the non-singular case
(see e.g. \cite{benedetti-petronio}).

\begin{cor}
Given a (topological) pants decomposition of $S$ with boundary 
curves $\gamma_1, \cdots, \gamma_N$, there is a homeomorphism 
$$ \cT_{S,x,\theta}\rightarrow (\R_{>0}\times \R)^N $$  
sending a hyperbolic metric to the length and fractional Dehn
twist parameters at the $\gamma_i$.
\end{cor}

The fractional Dehn twist parameters used here are the translation
length of one side with respect to the other so that, for a boundary
curve of length $l$, a parameter equal to $l$ corresponds to a
``usual'' Dehn twist (the other possibility
is to use an ``angle'' parameter, where $2\pi$ corresponds to 
full Dehn twist). 

\subsection{Proof of Proposition \ref{pr:qconf}}
\label{ap:A2}

It is now possible to use the pants decomposition provided by
Lemma \ref{lm:pant-existence} to prove Proposition \ref{pr:qconf}:
the induced metric on the boundary of the convex core is
(uniformly) quasi-conformal to the conformal structure at 
infinity. 

The starting point is that a pants decomposition of $(\dr M,m)$
with boundary curves of bounded length defines a pants decomposition
of $(\dr M, \tau)$ with boundary curves of approximately the same 
length. Recall that the constant $C_p$ was introduced in Lemma
\ref{lm:pant-existence}.

\begin{lemma} \label{lm:qc-pants}
There exists a constant $C>0$ as follows. 
Let $\gamma_1, \cdots, \gamma_N$ be simple closed curves on 
$\dr M$, defining a pants decomposition, of lengths less than 
$C_p$ for $m$. Then
$$ \forall i\in \{ 1,\cdots, N\}, \frac{L_m(\gamma_i)}{C}
\leq L_\tau(\gamma_i)\leq CL_m(\gamma_i)~. $$
\end{lemma}

\begin{proof}
The upper bound is a direct consequence of the first point in
Proposition \ref{pr:longueur-maj}. If $\gamma_i$ is short for $m$ --- i.e., it
is the core of a long tube in the thin part of $(\dr M,m)$ ---
then the second point of Proposition \ref{pr:longueur-maj} proves 
the lower bound for $\gamma_i$.

Suppose now that $\gamma_i$ is realized in $(\dr M, m)$ as a 
closed geodesic in the thick part of $\dr M$. Then there exists
a closed geodesic $\gamma'$ intersecting $\gamma_i$ of length
at most $C_p$. If the length of $\gamma_i$ in $(\dr M,\tau)$ were
small, than $\gamma_i$ would be realized in $(\dr M,\tau)$ as 
the core of a long tube $T$ in the thin part of $(\dr M,\tau)$. 
But then $\gamma'$ would have to be long (at least as long as 
the $T$). This would contradict the first point in Proposition 
\ref{pr:longueur-maj}, and this proves the lower bound for 
$\gamma_i$.
\end{proof}

\begin{lemma} \label{lm:twist}
There exists a constant $C>0$ such that, for each of the 
$\gamma_i$, the difference in the Dehn twist parameter 
corresponding to $\gamma_i$ in $m$ and in $\tau$ is
at most $C (|\log(L_m(\gamma_i))|+1)$. 
\end{lemma}

The precise form of the estimate is important only if $\gamma_i$
is short for $m$ (and therefore for $\tau$), in which case 
$|\log(L_m(\gamma_i))|$ is half the length of the tube in the thin
part of $(\dr M, m)$ containing $\gamma_i$.

\begin{proof}
Suppose first that $\gamma_i$ is not short. Then it is contained
in the thick part of $(\dr M, m)$, and there exists another 
curve $\gamma'$, intersecting $\gamma_i$, of uniformly bounded
length. A Dehn twist parameter bigger 
than some constant would extend the length of $\gamma'$ by 
more than is allowed by Proposition \ref{pr:longueur-maj}, 
this proves the lemma in this first case. 

The same argument can be used when $\gamma_i$ is short (i.e. when it
is the core of a long thin tube), then $\gamma'$ can be chosen to 
have length bounded by a constant time $|\log(L_m(\gamma_i))|$,
and this defines the maximal Dehn twist parameter along $\gamma_i$.
\end{proof}

\begin{proof}[Proof of Proposition \ref{pr:qconf}]
Let $\gamma_1, \cdots, \gamma_N$ be disjoint closed curves, defining
a pants decomposition of $(\dr M,m)$ with boundary curves of length
less than $C_p$, as in Lemma \ref{lm:pant-existence}. Let $l_1,
\cdots, l_N$ be the length of the $\gamma_i$ for $m$, and let
$d_1, \cdots, d_N$ be the Dehn twist parameters for the same
curves. 

Lemma \ref{lm:pant-topo} shows that the $\gamma_i$ also define a pant
decomposition of $(\dr M, \tau)$, let $l'_i$ be the length of the
$\gamma_i$ for $\tau$, and let $d'_i$ be their Dehn twist parameters.
Lemma \ref{lm:qc-pants} indicates that the $l'_i$ are within a fixed
multiplicative constant from the $l_i$, while, by Lemma 
\ref{lm:twist},
\begin{equation} \label{eq:dd}
|d'_i-d_i|\leq C (|\log(l_i)|+1)~,
\end{equation}
where $C$ is some positive constant.

Let $m'$ be the hyperbolic metric with cone singularities 
obtained by gluing pairs of pants with boundary lengths equal to 
the $l_i$, but with Dehn twist parameters equal to the $d'_i$. 

Note that $m'$ is $C_1$-quasi-conformal to $m$, for some 
uniform constant $C_1>0$. To prove this remark that for 
each $i\in \{ 1,\cdots, N\}$ the set of points at distance at
most $C(|\log(L_m(\gamma_i))|+c_M)$ from $\gamma_i$ is an annulus, and that those
annuli are disjoint. One can therefore build a $C_1$-quasi-conformal
diffeomorphism between $m$ and $m'$ which is an isometry in
the complement of those annuli around the $\gamma_i$, and which
is ``twisted'' in those annuli, with a twisting parameter which
is an affine function of the distance to the $\gamma_i$.

The second and last step is that 
$m'$ is $C_3$-quasi-conformal $\tau$. Since those two metrics differ
only by the lengths of the boundary curves $\gamma_i$, and in view
of (\ref{eq:dd}), this follows again from a simple and explicit 
construction which we leave to the interested reader.
\end{proof}

Note that it might be possible to prove Proposition \ref{pr:qconf} using
the same type of arguments as those used by Epstein and Marden \cite{epstein-marden} in the
non-singular case. This would have the advantage of providing directly a quasiconformal
constant independent on the genus of the boundar.